\documentclass[hidelinks,12pt,a4paper]{article} 

\usepackage[english]{babel}
\usepackage{latexsym}
\usepackage{amssymb}
\usepackage{amsmath}
\usepackage{amsfonts}
\usepackage{mathtools}
\usepackage{graphics}
\usepackage{enumitem}
\usepackage{afterpage}
\usepackage{fancybox}
\usepackage{hyperref}
\usepackage{ulem}
\usepackage{xspace}
\usepackage{xcolor}
\usepackage{color}
\usepackage{scalerel}
\usepackage{tikz}
\usetikzlibrary{graphs}
\usetikzlibrary{arrows}
\usetikzlibrary {quotes}

\newcommand{\al}{\alpha}
\newcommand{\be}{\beta}
\newcommand{\ga}{\gamma}
\newcommand{\de}{\delta}
\newcommand{\ep}{\varepsilon}
\newcommand{\la}{\lambda}
\newcommand{\fhi}{\varphi}

\newcommand{\ita}{\textit}

\newcommand{\mN}{\mathbb{N}}

\newcommand{\mT}{\mathbb{T}}

\newcommand{\sF}{\mathcal{F}}
\newcommand{\sG}{\mathcal{G}}

\newcommand{\sP}{\mathcal{P}}

\newcommand{\fX}{\mathfrak{X}}

\newcommand{\imp}{\Rightarrow}
\newcommand{\frec}{\rightarrow}

\newcommand{\x}{\cdot}
\newcommand{\vir}{``}
\newcommand{\sse}{\leftrightarrow}
\newcommand{\et}{\wedge}
\newcommand{\vel}{\vee}

\newcommand{\qed}{\hfill$\square$}

\newcommand{\g}{\ita{\texttt{g}}}
\newcommand{\wqo}{\textsf{wqo}\xspace}
\newcommand{\proof}{\noindent \ita{Proof:}\ }

\newcommand{\KTlw}{K$\mbox{T}_\ell(\omega)$}
\newcommand{\KTln}{K$\mbox{T}_\ell(n)$}
\newcommand{\KTw}{K$\mbox{T}(\omega)$}
\newcommand{\KTn}{K$\mbox{T}(n)$}

\newcommand{\RCA}{\textsf{RCA$_0$}}
\newcommand{\RCAs}{\textsf{RCA$^{*}_0$}}
\newcommand{\ACA}{\textsf{ACA$_0$}}

\newcommand{\ERA}{\ita{ERA}\ }

\renewcommand{\preceq}{\preccurlyeq}
\renewcommand{\leq}{\leqslant}
\renewcommand{\geq}{\geqslant}
\renewcommand{\underline}[1]{\textsc{#1}}
\newcommand{\defFoAs}{\:\leftrightarrow \! \! \raisebox{-0.5 ex}[0 ex][0 ex]{\tiny Def}\:}
\newcommand{\putAs}{\! \coloneqq}

\newtheorem{theor}{Theorem}[section]
\newtheorem{prop}{Proposition}[section]
\newtheorem{lem}{Lemma}[section]
\newtheorem{cor}{Corollary}[section]

\newtheorem{defi}{Definition}[section]
\newtheorem{conj}{Conjecture}[section]

\newtheorem{rem}{Remark}

\newcommand{\Rem}[1]{Remark~\ref{#1}}
\newcommand{\Theor}[1]{Theorem~\ref{#1}}
\newcommand{\Lem}[1]{Lemma~\ref{#1}}
\newcommand{\Def}[1]{Def.~\ref{#1}}

\newcommand{\Cor}[1]{Corollary~\ref{#1}}
\newcommand{\Prop}[1]{Proposition~\ref{#1}}
\newcommand{\Sec}[1]{Sec.~\ref{#1}}
\newcommand{\Subsec}[1]{Subsec.~\ref{#1}}

\newcommand{\COMMENT}[1]{}

\definecolor{navy}{rgb}{0,0,0.5}
\definecolor{magenta}{rgb}{0.7,0,0.7}
\graphicspath{{../grafici/}}

\title{Ordinal Analysis of Well-Ordering Principles, \\ Well Quasi-Orders Closure Properties, \\ and $\Sigma_n$-Collection Schema}
\author{Gabriele Buriola$^1$, Andreas Weiermann$^2$}
\date{%
    $^1$University of Verona\\ \vspace{0.1 cm}%
    $^2$Ghent University
}
\begin{document}

\maketitle

%\tableofcontents

\begin{abstract}

The study of well quasi-orders, {\wqo}, is a cornerstone of combinatorics and within {\wqo} theory Kruskal's theorem plays a crucial role. Extending previous proof-theoretic results, we calculate the $\Pi^1_1$ ordinals of two different versions of labelled Kruskal's theorem: $\forall n \,$\KTln \ and \KTlw \, denoting, respectively, all the cases of labelled Kruskal's theorem for trees with an upper bound on the branching degree, and the standard Kruskal's theorem for labelled trees.

In order to reach these computations, a key step is to move from Kruskal's theorem, which regards preservation of {\wqo}'s, to an equivalent Well-Ordering Principle (WOP), regarding instead preservation of well-orders. Given an ordinal function $\g$, WOP$(\g)$ amounts to the following principle $\forall \fX\, [\mbox{WO}(\fX) \frec \mbox{WO}(\g(\fX))]$, where WO$(\fX)$ states that ``$\fX$ is a well-order''. In our case, the two ordinal functions involved are $\g_{\forall}(\fX)=\sup_{n}\vartheta(\Omega^n \x \fX)$ and $\g_{\omega}(\fX)=\vartheta(\Omega^{\omega}\! \x \fX)$.

In addition to the ordinal analysis of Kruskal's theorem and its related WOP, a series of Well Quasi-orders Principles (WQP) is considered. Given a set operation $\sG$ that preserves the property of being a {\wqo}, its Well Quasi-orders Principle, WQP$(\sG)$, is given by the statement $\forall Q\, [Q\, \mbox{\wqo} \frec \sG(Q)\, \mbox{\wqo}]$. Conducting this study, unexpected connections with different principles arising from Ramsey and Computational theory, such as RT$^2_{<\infty}$, CAC, ADS, RT$^1_{<\infty}$, turn up.

Lastly, extending and combining previous results, we achieve also the ordinal analysis of the collection schema $\mbox{B}\Sigma_n$.

\medskip

\noindent{\bf MS Classification 2020:}\ \textbf{03B30}, 03F15, 03F35, 06A07.

\end{abstract}

\smallskip

%\footnotesize{bla}
%\noindent \textbf{MS Classification 2020:}\ NNNN, MMMM.
%\noindent \textbf{Keyword}: XXX, YYY, ZZZ.

\section{Introduction}\label{sec:introduction}

This paper deals with the ordinal and proof-theoretic analysis of a series of principles related to ordinals and well quasi-orders. More precisely, we treat: the well-ordering principles relative to four different ordinal functions, two versions of Kruskal's theorem for labelled trees, and the properties of well quasi-orders of being closed under, respectively, iterated disjoint union and product. On top of that, also the proof-theoretic ordinal of the collection schema $\mbox{B}\Sigma_n$ is computed.

\subsection{Content and Structure}\label{subsec:content}

Kruskal's tree theorem \cite{Kruskal60} is a milestone in the theory of well quasi-orders, {\wqo}, with ramified applications in many different areas, such as term rewriting \cite{Dershowitz82,Dershowitz87} and mathematical logic \cite{Gallier91,Simpson85nonprovability}. Extending previous proof-theoretic results on this topic (mainly \cite{RW93}), we calculate the $\Pi^1_1$ proof-theoretic ordinals of two different versions of labelled Kruskal's theorem. 

More precisely, the following two statements regarding labelled trees are considered. The former, denoted by \KTlw, is standard Kruskal's theorem stating that: ``for all $Q$ if $Q$ is a {\wqo}, then the set of finite trees with labels in $Q$, $\mT(Q)$, is a {\wqo}''. The latter, denoted by $\forall n \,$\KTln, concerns a parametrized version of the former. If $\mT_n(Q)$ is the set of finite trees labelled by $Q$ with branching degree at most $n$ (i.e., every node has at most $n$ children), we can consider the following restricted parametrized version of Kruskal's theorem \KTln: ``for all $Q$ if $Q$ is a {\wqo}, then $\mT_n(Q)$ is a {\wqo}''; then, $\forall n \,$\KTln \ is the union of all such restricted cases, namely $\{$\KTln $|\, n \! \in \! \mN \}$. For what concerns the corresponding ordinals, the following estimations are obtained: $|\mbox{RCA}_0 +\mbox{KT}_{\ell}(\omega)|=\vartheta (\Omega^{\omega+1})$ and $|\mbox{RCA}_0 + \forall n\, \mbox{KT}_{\ell}(n)|=\vartheta (\Omega^{\omega}+\omega)$. These two ordinals imply that the conjunction of all the finite cases, i.e., $\forall n\,$\KTln , is strictly weaker then the infinite case, \KTlw; thus, the situation differs from the unlabelled case, thoroughly treated by Michael Rathjen and the second author in their proof-theoretical investigations \cite{RW93}, where, over RCA$_0$,  $\forall n\,$KT$(n)$ and KT$(\omega)$ turn out to be equivalent. 

In order to reach the above calculations, a key step is to move from Kruskal's theorem, which regards preservation of {\wqo}'s, to an equivalent Well-Ordering Principle (WOP), regarding instead preservation of well-orders; for an introduction to WOP see \cite{AR09} and \cite{RW11}. Roughly speaking, given an ordinal function $\g$, WOP$(\g)$ amounts to the following principle $\forall \fX\, [\mbox{WO}(\fX) \frec \mbox{WO}(\g(\fX))]$, where WO$(\fX)$ is a formula stating that ``$\fX$ is a well-order''. In our case, the two ordinal functions involved are respectively $\g_{\omega}(\fX)=\vartheta(\Omega^{\omega}\! \x \fX)$ and $\g_{\forall}(\fX)=\sup_{n}\vartheta(\Omega^n \x \fX)$; the collapsing function $\vartheta$ is introduced in \Subsec{subsec:preliminaries} (for its definition, as well as other collapsing functions, see also \cite{Buchholz86,RW93,VanderMereen15}). The main tool used to achieve a proper analysis of the two related WOP's is the extension of a result, due to Arai \cite[Theorem 3]{Arai20}, regarding the proof-theoretic ordinal of a well-ordering principle. In particular, keeping the same thesis, we weak the hypotheses of the theorem to include the two ordinal functions risen from Kruskal's theorem. 

Beside the ordinal analysis of Kruskal's theorem and its related WOP, we study a series of Well Quasi-orders closure Properties (WQP), sometimes called Well Quasi-orders Principles; these can be considered the {\wqo} version of WOP's. More precisely, let $\sG$ be a set operation that preserves the property of being a {\wqo}, such as the Cartesian product, then the Well Quasi-orders closure Property associated to $\sG$, WQP$(\sG)$, is given by a suitable formalization of $\forall Q\, [Q\, \mbox{\wqo} \frec \sG(Q)\, \mbox{\wqo}]$. Conducting this study, unexpected connections with different principles arising from Ramsey and Computational theory, such as RT$^2_{<\infty}$, CAC, ADS, RT$^1_{<\infty}$, turn up; we briefly recall these principles in \Subsec{subsec:preliminaries}, but for a thorough introduction and analysis of their proof-theoretic relations we refer to \cite{HS07}.

Furthermore, extending and combining previous results, mainly \cite{Avigad:without} and \cite{Beklemishev:collection}, we achieve also the ordinal analysis of the collection schema $\mbox{B}\Sigma_n$ which, despite its relevance, seems to be missing in the literature. In particular, the following estimation is obtained $|\mbox{B}\Sigma_n|=\omega_n$, where $\omega_n$ is the $\omega$-tower of height $n$.

\smallskip

The structure of the article is as follows: the remaining part of this section exposes some general preliminaries; \Sec{sec:ordinalWOP} is devoted to extend Arai's aforementioned theorem to a wider class of ordinal functions and to apply the extended result to some specific well-ordering principles; in \Sec{sec:wqo}, which contains the main achievements of the paper, we compute the $\Pi^1_1$ proof-theoretic ordinals of the theories RCA$_0$ + \KTlw \ and RCA$_0$ + $\forall n\!$ \KTln, treating also the ordinal analysis of different well quasi-orders closure properties related to the product and the disjoint union of {\wqo}; \Sec{sec:Bsigma} is dedicated to the ordinal analysis of the collection schema $\mbox{B}\Sigma_n$ and regards mainly first-order logic; \Sec{sec:works} compares the contemporary literature and sketches future developments; finally, an ending section gives a comprehensive summary of the whole paper and draws conclusions.

Part of the material in this article has been treated in the PhD thesis of the first author \cite{BuriolaPhD}.

\subsection{Preliminaries}\label{subsec:preliminaries}

Before treating the first main topic, namely well-ordering principles, we settle down some concepts frequently recurring in the paper.

For what concerns notation, we adopt the standard one routinely used in mathematical logic, particularly in second-order arithmetic. For example, we use small Latin letters ($n,m,x,y,z, \dots$) for numeric variables, capital Latin letters ($X,Y,\dots$) for set variables, and Greek letters ($\al, \be, \ga, \dots , \omega, \Omega$) for ordinals. The hierarchies for $\Pi^0_n$, $\Sigma^0_n$, $\Pi^1_n$ or $\Sigma^1_n$ formulas are defined as usual; namely a formula $\fhi$ is in $\Sigma^0_0$, or equivalently in $\Pi^0_0$, if it has only bounded quantifiers and it is in $\Sigma^0_n$ (resp. in $\Pi^0_n$) if it is of the form $\exists x\, \theta(x)$ (resp. $\forall x\, \theta(x)$) with $\theta(x)$ in $\Pi^0_{n-1}$ (resp. $\Sigma^0_{n-1}$). The hierarchies $\Pi^1_n$ and $\Sigma^1_n$ are defined similarly, considering this time set quantifiers; in \Sec{sec:Bsigma} the corresponding hierarchies in first-order logic are used. Moreover, we denote with $\mN$ the set of natural number as seen from within a second-order theory, i.e., the set $\{ x\, | \, x=x \}$ which exists by $\Sigma^0_0$-comprehension, and with $\omega$ the ``true'' meta-theoretic set of natural numbers.

\subsubsection{Reverse Mathematics}\label{subsub:reverse}

The \ita{Reverse Mathematics} programme, pioneered by Harvey Friedman \cite{Friedman75}, and subsequently developed by Stephen Simpson and others \cite{Simpson05,Simpson09,Stillwell18}, aims to classify ``ordinary mathematics'' statements using as benchmark suitable axioms, mainly existential axioms, in the language of second-order arithmetic. Given the numerous fields of mathematics, the overall picture might seem to lack any regularity; however, most theorems in ordinary mathematics turn out to be provable in the base theory \RCA \, or to be equivalent, over \RCA, to one of four specific subsystems. Moreover, these systems, called the ``Big Five'' of reverse mathematics (in increasing order of strength: \textsf{RCA$_0$}, \textsf{WKL$_0$}, \textsf{ACA$_0$}, \textsf{ATR$_0$} and \textsf{$\Pi^1_1$-CA$_0$}), turn out to be linearly ordered. For a comprehensive introduction we refer to \cite{Simpson09}.

Here, we limit ourself to shortly expose the three subsystems of second-order arithmetic that are relevant for our investigations: \RCA, \RCAs, and~\ACA.

\smallskip

\noindent \textbf{\textsf{RCA$_0$}} is the weakest of the ``big five'' subsystems of reverse mathematics and is routinely used as base theory to compare stronger systems or results. Besides the basic axioms for natural numbers (e.g. $n+0=n, n\x 0=0, m +n =n+m,\dots$), its axioms contain the schema for $\Sigma^0_1$ induction
\[
\fhi(0) \, \et \, \forall n \left( \fhi(n) \frec \fhi(n+1) \right) \  \frec \ \forall n\, \fhi(n)
\]
for all $\Sigma^0_1$ formulas $\fhi$, and the schema for recursive, or $\Delta^0_1$, comprehension
\[
 \forall n \left( \fhi(n) \leftrightarrow \psi(n) \right) \ \frec \ \exists X\, \forall n\, \left( n \! \in \! X \leftrightarrow \fhi(n) \right)
\]
for all $\Sigma^0_1$ formulas $\fhi$ and all $\Pi^0_1$ formulas $\psi$.

Given its strong connections with computability theory and Turing machines, and despite the use of excluded middle and other differences \cite[Remark I.8.9]{Simpson09}, RCA$_0$ corresponds to some extend to Bishop's constructive mathematics \cite{Bishop67}(regarding this correspondence see also \cite{FSS83}). Finally, for what concerns an ordinal measure of the strength of RCA$_0$, the $\Pi^1_1$ ordinal\footnote{See \Def{def:Pi11Ord}.} of RCA$_0$ is $\omega^{\omega}$.

\smallskip

\noindent \textbf{\textsf{RCA$^*_0$}} is a subsystem weaker than \RCA\, and is commonly used to study results that are provable in \RCA. Differently from the other subsystems of second-order arithmetic, the language $L^2_{exp}$ of \RCAs \, includes a binary function $exp(n,m)$ representing exponentiation. Besides the basic axioms for natural numbers and the ones for exponentiation, namely $exp(n,0)=1$ and $exp(m,n+1)=exp(m,n)\x m$, \RCAs \, includes the induction schema for $\Sigma^0_0$ formulas of $L^2_{exp}$ and the $\Delta^0_1$ comprehension schema for $L^2_{exp}$ formulas. Its $\Pi^1_1$ ordinal is $\omega^3$.

\smallskip

\noindent \textbf{\textsf{ACA$_0$}}, which stands for Arithmetical Comprehension Axiom, is the third subsystem of the Big Five and, together with \RCA, plays a prominent role in reverse mathematics. Besides the basic ones for natural numbers, its axioms are given by the induction axiom
\[
0 \! \in \! X \, \et \, \forall n \left( n \! \in \! X \, \frec \, n+1 \! \in \! X \right) \ \frec \ \forall n\, n \! \in \! X
\]
together with the comprehension schema $\exists X\, \forall n\, (n\! \in \! X \leftrightarrow \phi(n))$ for all arithmetical formulas, i.e., formulas without set quantifiers.

\textsf{ACA$_0$} is tightly connected with Peano Arithmetic \cite[Remark I.3.3]{Simpson09} (see also \cite{HP98}). More precisely, \textsf{ACA$_0$} is a conservative extension of PA; this means that, for any sentence $\sigma$ in the language of first-order arithmetic, $\sigma$ is a theorem of PA if and only if $\sigma$ is a theorem of \textsf{ACA$_0$}. Said in other words, PA is the first-order part of \textsf{ACA$_0$}. In particular, they share the same $\Pi^1_1$ ordinals, namely $\ep_0$.

\subsubsection{Well-orders and the $\vartheta$ collapsing function}\label{subsub:wo}

Treating well-orders in the framework of reverse mathematics, we briefly recall how to encode them using second-order arithmetic.

Given a binary relation $\prec$ and a subset $Y$, we define:
\[
\begin{array}{rclcl} \vspace{0.04 cm}
Prog(\prec,Y) & \putAs & \forall x\, \left[ \forall y\, (y \prec x \frec y \! \in \! Y) \frec x \! \in \! Y \right] & & (\mbox{\ita{progressiveness}}) \\ \vspace{0.04 cm}
TI(\prec,Y) & \putAs & Prog(\prec,Y) \frec \forall x\, x \! \in \! Y & & (\mbox{\ita{transfinite induction}}) \\ \vspace{0.04 cm}
WF(\prec) & \putAs & \forall Y\, TI(\prec,Y) & & (\mbox{\ita{well-foundedness}}) \\
WO(\prec) & \putAs & WF(\prec) \, \et \, LN(\prec) & & (\mbox{\ita{well-orderedness}})
\end{array}
\]
where $LN(\prec)$ is a sentence asserting that $\prec$ is a linear order. Moreover, since any subset $X \subseteq \mN$ gives rise to a binary relation putting $x \leq_X y \defFoAs \langle x,y \rangle \in X$, in the following we use also the notation WO$(X)$.

\begin{rem}\label{rem:TIfirst} An appropriate definition of $TI(\prec,U)$ can be given also in first-order logic using an uninterpreted unary relation symbol $U;$ see \cite{Rathjen99}.

\end{rem}

Given its role in our results, let us briefly recall the definition, together with the main properties, of the collapsing function $\vartheta$; this part is mainly taken from \cite{RW93} to which we refer for more details as well as for the corresponding ordinal notation.

Let $E$ be the set of $\ep$-numbers (i.e., ordinals $\al$ such that $\omega^{\al}=\al$), and $\Omega$ be the first uncountable ordinal; for any $\al < \ep_{\Omega+1}$, we recursively define the set $E_{\Omega}(\al)$ whose elements are the $\ep$-numbers below $\Omega$ needed in the Cantor normal form of $\al$, denoted here by $=_{NF}$. $E_{\Omega}(\al)$ is given by:
\begin{enumerate}

\item $E_{\Omega}(0) \putAs E_{\Omega}(\Omega) \putAs \emptyset$,

\item $E_{\Omega}(\al) \putAs \al$, if $\al \in E \cap \Omega$,

\item $E_{\Omega}(\al) \putAs E_{\Omega}(\be) \cup E_{\Omega}(\delta)$ if $\al =_{NF} \omega^{\be} + \delta$.

\end{enumerate}
Let $\al^* \putAs \max\{E_{\Omega}(\al) \cup \{0\} \}$.

Next, given two ordinals $\al, \be < \ep_{\Omega+1}$, we define sets of ordinals $C_n(\al, \be)$, $C(\al, \be)$ and the ordinal $\vartheta(\al)$ by main recursion on $\al$, and side recursion on $n < \omega$, as follows:
\begin{itemize}

\item[C1)] $\{0,\Omega\} \cup \be \subseteq C_n(\al, \be)$ for all $n < \omega$,

\item[C2)] $\ga, \delta \in C_n(\al,\be) \, \et \, \xi =_{NF} \omega^{\ga} + \de \ \imp \ \xi \in C_{n+1}(\al, \be)$,

\item[C3)] $\de \in C_n(\al,\be) \cap \al \ \imp \ \vartheta(\de) \in C_{n+1}(\al,\be)$,

\item[C4)] $C(\al, \be) \putAs \bigcup \{ C_n(\al,\be) \, | \, n < \omega\}$,

\item[C5)] $\vartheta(\al) \putAs \min\{\xi < \Omega \, | \, C(\al, \xi) \cap \Omega \subseteq \xi \, \et \, \al \in C(\al, \xi) \}$.

\end{itemize}

Let us summarize the main properties of $\vartheta(\al)$, see \cite[Lemma 1.1 and Lemma 1.2]{RW93} for the proofs.

\begin{lem}\label{lem:vartheta} Let $\vartheta(\al)$ be defined as before, then:
\begin{enumerate}[label=(\arabic*)]

\item $\vartheta(\al)$ is defined for every $\al < \ep_{\Omega+1}$,

\item $\vartheta(\al) \in E$,

\item $\al \in C(\al, \vartheta(\al))$,

\item $\vartheta(\al) = C(\al, \vartheta(\al)) \cap \Omega$, and $\vartheta(\al) \notin C(\al, \vartheta(\al))$,

\item $\ga \in C(\al, \be) \, \Leftrightarrow \, E_{\Omega}(\ga) \subseteq C(\al, \be)$,

\item $\al^* < \vartheta(\al)$,

\item $\vartheta(\al) = \vartheta(\be) \, \Leftrightarrow \, \al = \be$,

\item $\vartheta(\al) < \vartheta(\be) \, \Leftrightarrow \, (\al \! < \! \be \et \al^* < \vartheta(\be)) \, \vel \, (\be \! < \! \al \et \vartheta(\al) \leq \be^*)$,

\item $\be < \vartheta(\al) \, \Leftrightarrow \, \omega^{\be} < \vartheta(\al)$. 

\end{enumerate}
\end{lem}

In the proof of \Prop{prop:g'0} the following, yet equivalent, definition of $\vartheta(\al)$ is exploited:
\[
\vartheta(\al) \! = \! \min \{\xi \! \in \! E \, |\, {\al}^* \! < \! \xi  \ \mbox{and}\ \forall \be \! < \! \al \, ( \be^* \! < \! \xi \frec \vartheta(\be) \! < \! \xi) \},
\]
see \cite[pag. 15]{VanderMereen15} for proof of the equivalence.

\subsubsection{Well Quasi-orders, Trees and Kruskal's Theorem}\label{subsub:Kruskal}

Some standard definitions regarding well quasi-orders and trees are here recalled, followed by the formal statement of Kruskal's theorem. For a recent introductory survey to the theory of well quasi-orders and better quasi-orders see \cite{BS:wqo}.

\begin{defi}\label{def:qo} A \textsc{quasi-order}, qo, $(Q, \leq)$ is given by a subset $Q \subseteq \mN$ together with a binary relation $\leq \, \subseteq Q \! \times \! Q$ which is reflexive and transitive; if in addition $\leq$ is antisymmetric, then $(Q,\leq)$ is a \textsc{partial order}, po.
\end{defi}
In the following,  we may denote with $Q$ also the qo $(Q,\leq)$, omitting the quasi-order relation $\leq$.

For well quasi-orders, the following definitions are adopted. 
\begin{defi}\label{def:wqo} Given a qo $(Q,\leq)$: 

\begin{enumerate}

 \item a sequence $(q_k)_k$ in $Q$ is \textsc{good} if there exist indexes $i <j$ such that $q_i \leq q_j$, otherwise is \textsc{bad};

 \item $(Q,\leq)$ is a \textsc{well quasi-order}, {\wqo}, if every infinite sequence in $Q$ is good, i.e., there are no infinite bad sequences.

\end{enumerate}

\end{defi}
We stress that there are actually many different definitions for {\wqo}, and that they are not equivalent with each other over weak theories; see \cite{CMS04,Marcone05,Marcone20} for a thorough analysis.

Let us briefly introduce finite ordered labelled trees.
\begin{defi}\label{def:trees} Given a set $Q$, the set $\mT(Q)$ of (finite ordered) \textsc{trees} with \textsc{labels} in $Q$ is recursively defined as follows:
\begin{enumerate}

\item for each $q \! \in \! Q$, $q[]$ is an element of $\mT(Q)$;

\item if $t_1, \dots , t_n$, with $n \! \geq \! 1$, is a finite sequence of elements of $\mT(Q)$ and $q \! \in \! Q$, then 
$t \putAs q[t_1, \dots , t_n]$ is an element of $\mT(Q)$.

\end{enumerate}
\end{defi}
Since all trees considered in this paper are finite and ordered, we omit these specifications.

Connected to labelled trees, there are the following definitions.
\begin{defi}\label{def:labels} Let $Q$ and $\mT(Q)$ be as before:

\begin{itemize}

\item if $t = q[t_1, \dots , t_n]$, then $q$ is the \textsc{label} of the root and $t_1, \dots , t_n$ are called the \textsc{(immediate) subtrees} of $t$;

\item if $Q$ is a singleton, then $\mT \putAs \mT(Q)$ is the set of \textsc{unlabelled} trees.

\end{itemize}
\end{defi}

Starting from a qo $Q$, it is possible to define an embeddability relation.
\begin{defi}\label{def:treemb} Given a qo $(Q,\leq_Q)$ and $t,s \! \in \! \mT(Q)$, we inductively define the embeddability relation $\preceq$ on $\mT(Q)$ as follows: $t \preceq s$ holds if

\begin{enumerate}

\item $t=p[]$, $s=q[]$ and $ p \leq_Q q$; or

\item $s=q[s_1, \dots , s_m]$ and $t \preceq s_i$ for some $ 1 \! \leq \! i \! \leq \! m$; or

\item $t=p[t_1, \dots, t_n]$, $s=q[s_1, \dots , s_m]$, $p \leq_Q q$ and there exist \\ $1  \leq \! i_1  < \dots <  i_n  \leq  m $ such that $t_k \preceq s_{i_k}$ for all $1  \leq  k  \leq  n$.

\end{enumerate}

\end{defi}
Let us observe that $(\mT(Q), \preceq)$ is a quasi-order.

Before stating Kruskal's theorems, one last definition is required: branching degree.
\begin{defi} The \textsc{branching degree} of a tree $t \! \in \! \mT(Q)$ is recursively defined as:
\begin{enumerate}

\item if $t=q[]$, then $Deg(t) \putAs 0$;

\item if $t=q[t_1, \dots, t_n]$, then $Deg(t) \putAs \max\{n, Deg(t_1), \dots, Deg(t_n)\}$.

\end{enumerate}
Moreover, for every $n \! \in \! \mN$ and every set $Q$:

\begin{itemize}

\item $\mT_n$ is the set of unlabelled trees with branching degree less or equal to~$n$;

\item $\mT_n(Q)$ is the set of trees with labels in $Q$ and branching degree less or equal to $n$.

\end{itemize}

\end{defi}
Since only finite trees are considered, $Deg$ is always well defined.

The various versions of Kruskal's theorem under scrutiny are now presented; each is denoted by an abbreviation, possibly accompanied by a numeric parameter.

\begin{theor}\label{theor:Kruskal} The following statements hold for every natural number $n$:

\begin{enumerate}

\item \KTn: $(\mT_n, \preceq)$ is a {\wqo};

\item \KTw: $(\mT,\preceq)$ is a {\wqo};

\item \KTln: if $Q$ is a {\wqo}, then $(\mT_n(Q), \preceq)$ is a {\wqo};

\item \KTlw: if $Q$ is a {\wqo}, then $(\mT(Q), \preceq)$ is a {\wqo}.

\end{enumerate}

\end{theor}

\subsubsection{Combinatorial principles from Computational and Ramsey theory}\label{subsub:combinatorial}

As we will see, although the key elements of this paper are well-ordering principles and well quasi-orders, unexpected connections with results stemming from computability and Ramsey theory arise. To allow also the reader not acquainted with such topics to take advantage of our results, the main principles involved in the article are collected below. These principles are part of the so-called ``zoo'' of reverse mathematics \cite{DM:logic, Freund:zoo}; we refer to \cite{HS07, SY:pairs} and the references therein for further readings.

For what concerns notation, given a set $A \subseteq \mN$ and $k\! \in \! \mN$, $[A]^k$ and $[A]^\omega$ denote respectively the subsets of $A$ with exactly $k$ elements and the infinite subsets of $A$; moreover, $k$ is identified with the set $\{0, \dots, k-1\}$.

The principles relevant for this paper are the following:

\smallskip

\noindent \textbf{RT$^2_{<\infty}$} Ramsey's Theorem for pairs and arbitrarily many colours: for every $k \! \in \! \mN$ and every $k$-colouring $c: [\mN]^2 \frec k$ of pairs there exists a $c$ homogeneous subset $H$; i.e., $c$ is constant on $[H]^2$.

%\smallskip

%\noindent \textbf{RT$^2_2$} Ramsey's Theorem for pairs and $2$ colours: every $2$-colouring $c: [\mN]^2 \frec \{0,1\}$ of pairs has homogeneous subset $H$.

\smallskip

\noindent \textbf{CAC} Chain Anti-Chain principle: every infinite partial order has either an infinite chain or an infinite antichain.

\smallskip

\noindent \textbf{ADS} Ascending Descending Sequence principle: every infinite linear order has either an infinite ascending sequence or an infinite descending sequence.

\smallskip

\noindent \textbf{COH} Cohesive principle: every sequence $\vec{R}=\{R_i \, | \, i \! \in \! \mN \}$ of subsets of $\mN$ has a $\vec{R}$-cohesive set; namely there is an infinite set $S$ such that
\[
\forall i\, \exists s\, \left[ \forall j \! > \! s\, (j \in S \frec j \in R_i) \, \vel \, \forall j \! > \! s\, (j \in S \frec j \notin R_i) \right].
\]

\smallskip

\noindent \textbf{RT$^1_{<\infty}$} Infinite Pigeonhole principle: for every $k\! \in \! \mN$ and every function $f\! : \mN \frec k$, there exist $n \! \in \{0, \dots k-1\}$ and an infinite set $A \subseteq \mN$ such that $f(i)=n$ for all $i\! \in \! A$.

For what concerns their proof-theoretic relations over the base theory \RCA, all the possible implications are exhaustively exposed in the figure below (see \cite{HS07, SY:pairs} for the proofs); in particular no arrow can be inverted.

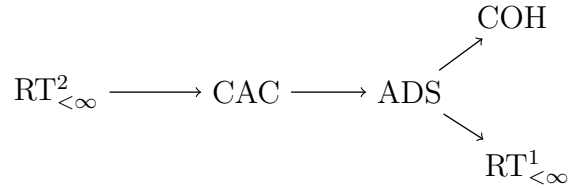
\begin{figure}[!htb]
\begin{center}
\begin{tikzpicture}[every node/.style={rectangle, fill=white, text=black}]
\node (RT2) at (-5,0) {RT$^2_{<\infty}$};
%\node (RT22) at (-4.25,0) {RT$^2_2$};
\node (CAC) at (-2.5,0) {CAC};
\node (ADS)    at (-0.35,0)  {ADS};
\node (COH)   at (1,1)  {COH};
\node (RT1)    at (1.2,-1)  {RT$^1_{<\infty}$};

%\draw[line width=0.5pt,  ->] (RT2)--(RT22);
%\draw[line width=0.5pt,  ->] (RT22)--(CAC);
\draw[line width=0.5pt,  ->] (RT2)--(CAC);
\draw[line width=0.5pt,  ->] (CAC)--(ADS);
\draw[line width=0.5pt,  ->] (ADS)--(COH);
\draw[line width=0.5pt,  ->] (ADS)--(RT1);

\end{tikzpicture}

\end{center}

\end{figure}

\section{Ordinal Analysis of Well-Ordering Principles}\label{sec:ordinalWOP}

We briefly recall two central concepts in our investigation, \ita{Ordinal Analysis} and \ita{Well-Ordering Principles}.

The origin of ordinal analysis can be tracked back to Gentzen \cite{Gentzen36, Gentzen43} who showed that transfinite induction up to the ordinal
\[
\ep_0 = \sup \{ \omega, \omega^{\omega} \! , \omega^{\omega^{\omega}}\! , \dots \}
\]
suffices to prove the consistency of PA. Moreover, he proved that $\ep_0$ is the best possible choice in the sense that PA proves transfinite induction up to $\al$ for arithmetic formulas for any $\al < \ep_0$. Thus, the idea is that $\ep_0$ in some sense ``measures'' the consistency strength of PA.

From the seminal works of Gentzen, ordinal analysis has grown enormously both in results and methods; still, one of the goals, perhaps the best known, is roughly to ``attach ordinals in a given representation system to formal theories'' \cite[Pag. 1]{Rathjen99}. One of the most commonly attached ordinals is the $\Pi^1_1$ ordinal; following \cite[Definition 1.11]{Arai20book}, we define it as follows:
%which, for many ``natural'' theories $T$,\footnote{A formal and rigorous definition of naturalness for a logical theory is still an open problem.} equals the supremum of the provable recursive well-orderings of $T$, namely
\begin{defi}\label{def:Pi11Ord} For a theory T, %comprising $RCA^*_0$, 
the $\Pi^1_1$ ordinal $|T|_{\Pi^1_1}$ is defined as:
\[
|T|_{\Pi^1_1} \putAs \sup \{ ot(\prec) \mid \, \prec \, \mbox{is recursive and}\ T \vdash \mbox{WO}(\prec) \},
\]
where $ot(\prec)$ is the \ita{order type} of $\prec$ (see \cite{dJP77,Schmidt79}), and WO$(\prec)$ is the previously defined sentence asserting that $\prec$ is a well-order. 
\end{defi}
By Spector's $\Sigma^1_1$-boundedness theorem\footnote{See for example \cite[Theorem 1.5]{Arai20book} or \cite[Theorem 2.4]{Walsh:incompleteness}.}, if $T$ is a recursive $\Pi^1_1$-sound theory, then $\left|T\right|_{\Pi^1_1} \! < \! \omega^{ck}_1$; where the Church-Kleene ordinal $\omega^{ck}_1$ denotes the least non-recursive ordinal. For sake of readability, $\left| T \right|_{\Pi^1_1}$ is simply denoted by $|T|$; for further readings on ordinal analysis, we refer to \cite{Arai20book,Rathjen99,Rathjen06}.

Let us observe that, given \Rem{rem:TIfirst}, a suitable definition of $\left|T\right|_{\Pi^1_1}$ can be given also for a first-order theory $T$.

\smallskip

Well-ordering principles, which concern the preservation of well-orderedness under ordinal functions, are both a central topic of study and a crucial tool in our ordinal analysis of Kruskal's theorem. Their precise definition is as follows.

\begin{defi}\label{defi:WOP} Given an ordinal function $\g \!: \Omega \frec \Omega$, the \textsc{well-ordering principle} of $\g$, WOP(\g), amounts to the following statement
\begin{equation}\label{eq:WOP}
 \forall \fX \, \left[ \mbox{WO}(\fX) \frec \mbox{WO}(\g(\fX)) \right].
\end{equation}
\end{defi}
In order to study well-ordering principles in the framework of reverse mathematics, all ordinal functions $\g$ that give rise to a WOP considered in this paper are assumed to have a corresponding computable term structure; namely, all elements of $\g(\fX)$ are computably constructed starting from the ones of $\fX$ together with a finite list of computable functions. See \cite{Rathjen06} for more details regarding term structures for ordinal representations, and \Subsec{subsec:WOP} for the term structure of one of the ordinal functions under scrutiny.

As already mentioned in the introduction, a key step in our ordinal analysis of labelled Kruskal's theorem is to move from Kruskal's result, which is about preservation of {\wqo}, to an equivalent well-ordering principle, concerning instead well-orders. Although the literature regarding well-ordering principles is already well established \cite{AR09,Arai17,Arai20,Freund19,Girard87,MM11,Rathjen14,Rathjen22,RT20,RV15,RW11}, the WOP's required for our case, and their proof-theoretic analysis, have not been treated yet. The gap is filled in the next section where the following ordinal estimation
\[
|\mbox{\textsf{ACA}}_0 + \mbox{WOP}(\g)|=\g'(0) \putAs \sup_n \g^n(0)=\min\{\al \! > \! 0 \, | \, \forall \be \! < \! \al\ \g(\be) \! < \! \al \}
\]
is achieved for a class of ordinal functions $\g$ larger than the one considered by Arai \cite{Arai20}, an example of which is given in \Def{def:ordinalfun}.% We obtain the aforementioned estimation by extending a previous result due to Arai \cite[Theorem 3]{Arai20}.
 %moreover, preparing this section, we glimpsed another possible approach which instead uses an equivalence lemma, proved by Pakhomov and Walsh \cite[Lemma 3.8]{PW23}, relating well-ordering principles and well-ordering rules.

\subsection{Extending Arai's theorem}\label{subsec:arai}

In \cite{Arai20}, Arai studied the proof-theoretic ordinal of \ACA \, + WOP(\g) for a \ita{normal function} \g, namely for a strictly increasing and continuous ordinal function. More precisely, given a normal function $\g$ and denoting with $\g'$ its derivative (i.e., the ordinal function enumerating the fixed points of \g), if  \g \ and $\g'$ satisfy some term properties \cite[Def. 3 and Def. 4]{Arai20}, then the following ordinal computation holds \cite[Theorem 3]{Arai20}:
\[
|\mbox{\textsf{ACA}}_0 + \mbox{WOP}(\g)|={\g}'(0)=\min\{ \al  \, |\, \g(\al)=\al \}.
\]
We aim to extend this result, weakening the conditions required for \g.

\smallskip

Keeping the same notation as in \cite{Arai20}, we consider a computable function \g$: \sP(\mN) \frec \sP(\mN)$ which sends a linear order $<_{X}=\{(n,m) \,| \langle n,m \rangle \! \in \! X\, \}$ to a linear order $<_{\g(X)}=\{(n,m) \,| \langle n,m \rangle \! \in \! \g(X)\, \}$\footnote{By $\langle n,m \rangle$ we denote the encoding in $\mN$ of pairs of natural numbers. Moreover, with a slight and harmless abuse of notation, $\g$ indicates both the function from $\sP(\mN)$ to $\sP(\mN)$ and its restriction \g$: \Omega \frec \Omega$ which sends countable ordinals to countable ordinals.}. However, since we treat functions which may not be normal, ${\g}'(0)$ is defined not as $g$'s first fixed point (which might not exist), but directly as the first ordinal closed under \g, namely ${\g}'(0):=\sup_{n}{\g}^n(0)=\min\{ \al \! > \! 0 \, |\, \forall \be \! < \! \al\ \g(\be) \! < \! \al \}$. Differently from \cite{Arai20}, for $\g$ only the following properties are required:

\begin{enumerate}

\item \g\ is weakly increasing, i.e., $\al \! \leq \! \be \ \imp \ \g(\al) \! \leq \! \g(\be)$;

\item ${\g}'(0)$ is an epsilon number, i.e., $\omega^{{\g}'(0)}={\g}'(0)$.

\end{enumerate}
Thus, with respect to Arai's conditions for \g, we drop normality as well as the term structure requirements (for this point see also \Rem{rem:NotTerm}); let us emphasise how property 2. suffices to ensure that $\g'(0) \! > \! 0$.

\bigskip

Our extended version of \cite[Theorem 3]{Arai20} can now be stated:

\begin{theor}\label{theor:WOP} Given an ordinal function \g\ as above, then:
\[
|\mbox{\textsf{ACA}}_0 + \mbox{WOP}(\g)|={\g}'(0).
\]
\end{theor}

\noindent \ita{Proof:} The easy direction can be immediately proven considering the following ordinal succession: $\al_0 \putAs 0, \al_1\putAs \g(0)$, and $\al_{n+1} \putAs \g(\al_{n})$. By definition, we have that $\sup_{n}\al_{n}=\sup_{n}\g^n(0)=\g'(0)$; but, using finitely many iterations of WOP(\g), each $\al_n$ is a well-order and thus $|\mbox{\textsf{ACA}}_0 + \mbox{WOP}(\g)| \! \geq \g'(0)$.

\smallskip

For the other direction, i.e., $\leq$, we resort to a thorough analysis of \cite[Theorem 3]{Arai20} and its proof in order to extract the key points where the properties of \g\ are actually used, and subsequently, check that our weaker hypotheses for \g\ are indeed sufficient. Given the structure of our proof, all the references in the remaining part of this subsection hint at Arai's article \cite{Arai20}.

For sake of clarity, we briefly summarize Arai's proof (cf. \cite[pag. 266-268]{Arai20}). The idea is to move through three different calculi, $G_2 + (VJ) + (prg) + (WPL)$, $(prg)^{\infty}+(WP)$, and $Diag(\emptyset)+(prg)^{\emptyset}$, in order to extract, from the proof of the well-foundedness of $\prec$, an upper bound on the order type of $\prec$. Although different, the three calculi share some common features. In particular, $(prg), (prg)^{\infty},$ and $(prg)^{\emptyset}$ are progressive rules that turn progressiviness into an appropriate rule version; for example:
\[
(prg): \ \ \ \frac{\Gamma,\, E_A(t),\, x \, {\nless}_A \, t,\, E_A(x)}{\Gamma,\, E_A(t)}.
\]
$(WPL)$ and $(WP)$ are instead well-ordering rules for $\g$, namely the well-ordering principle WOP$(\g)$ in form of rule; e.g.
\[
(WPL): \ \ \ \frac{\Gamma, E_A(x)\  \ \ \ \Gamma,LO(<_A)\  \ \ \ \Gamma,\neg TI[<_{{\g}_A}] }{\Gamma} .
\]
Lastly, $G_2$ is a calculus for second-order logic; while $(VJ)$ and $Diag(\emptyset)$ are, respectively, an inference rule for the induction schema and the set of initial sequents of the calculus $Diag(\emptyset)+(prg)^{\emptyset}$. We refer to Arai's paper \cite{Arai20} for more details.

The proof proceeds as follows. Assume that \ACA\, + WOP(\g) proves WO($\prec$) for a linear relation $\prec$. One can obtain a derivation of $\Delta_0, E_{\prec}(x)$ in the calculus $G_2 + (VJ) + (prg) + (WPL)$, where $\Delta_0$ is a set of negated axioms and $E_{\prec}$ a fresh new variable related to $\prec$. Next, $G_2 + (VJ) + (prg) + (WPL)$ is embedded into $(prg)^{\infty}+(WP)+(cut)_{1^{st}}$, an intermediate infinitary calculus obtained from $(prg)^{\infty}+(WP)$ adding a first-order cut rule; namely, from $G_2 + (VJ) + (prg) + (WPL) \vdash \Delta_0, E_{\prec}(x)$, we move to $(prg)^{\infty}+(WP)+(cut)_{1^{st}}\vdash^{\omega^2}_{d,p} \Delta_0, E_{\prec}(n)$ for all $n$, where $\omega^2$ bounds derivation length, $d$ the number of nested applications of (WP) and $p$ the rank of cut formulas. 
Applying cut elimination, we arrive at $(prg)^{\infty}+(WP) \vdash^{\be}_{c} E_{\prec}(n)$ for all $n$, with $\be \! < \! \ep_0$ and $c\! < \! \omega$. Theorem 5 (which is a $(WP)$ elimination result) allows to move to $Diag(\emptyset)+(prg)^{\emptyset}$ in which $\{n\}\vdash^{\al} E_{\prec}(n)$ holds for all $n$, with $\al =F(\be,c)+\be$ and $F$ a suitable function. Finally, thanks to Theorem 6 and Proposition 2, we can extract an order-preserving injection $f \! : | \! \! \prec \! \! | \frec \omega^{\al+1}$ such that $n \prec m$ implies $f(n) < f(m) < \omega^{\al+1} = \omega^{F(\be,c)+\be+1} < \omega^{\g'(0)}=\g'(0)$; thus $ot(\prec) \! \leq \! \g'(0)$.

All in all, Arai's proof revolves around three main ingredients: three different calculi for second-order arithmetic, two of which infinitary; a bounding function $F$ defined from \g , used to control derivation length and whose required properties are ensured by Arai's Proposition 2 [pag. $265$], and two theorems, Theorem 5 and Theorem 6 [pag. $265$ and pag. $266$] which are used to move from one calculus to another and to extract, together with $F$, the final upper bound for the proof-theoretic ordinal. We consider these elements one by one.

For what concerns the three main calculi used in the proof\,\footnote{In Arai's article, the last two calculi, which are the infinitary ones, refer also to a subsidiary set $\sP$; this set is needed only for the second main result of the paper (\cite[Theorem 4]{Arai20}) and thus it is negligible here.}, $G_2 + (VJ) + (prg) + (WPL)$, $(prg)^{\infty}+(WP)$ and $Diag(\emptyset) + (prg)^{\emptyset}$, these are defined without any concrete reference to the properties of \g; thus, we can keep them as they are.

For the function $F$ Arai's definition, which we briefly recall, is still appropriate. Given \g, we define $\be, \al \mapsto F(\be, \al)$ by induction on $\al$ as follows: $F(\be,0)=\omega^{1+\be}$, $F(\be, \al+1)=F(\g(\omega^{F(\be, \al) +\be +1}), \al)+1$, and $F(\be, \la)=\sup\{F(\be,\al)+1 \, |\, \al \! < \! \la\}$ for $\la$ limit ordinal. From our hypotheses for \g, the following weaker, yet sufficient, version of Arai's Proposition 2 can be obtained:

\noindent \textbf{Claim.} If $F(\be, \al)$ is defined as above, then:

\begin{enumerate}

\item if $\al \! \leq \! \ga$ and $\be \! \leq \! \de$, then $F(\be, \al) \! \leq \! F(\de, \ga)$;%, i.e., $F$ is weakly increasing in each argument;

\item if $\be \! < \! \g'(0)$ and $c \! < \! \omega$, then $F(\be, c) \! < \! \g'(0)$.

\end{enumerate}

\noindent \ita{Proof of the Claim.} 1. derives from the fact that each function $\be \mapsto \be + \al$, $\be \mapsto \omega^{\be}$ and $\be \mapsto \g(\be)$ is weakly increasing; 2. is proved by induction on $c$, using the closure of $\g'(0)$ with respect to \g. \ita{This concludes the proof of the Claim.}

Regarding the role of Proposition 2 in Arai's proof, respectively the previous Claim in ours, 1. is extensively used to majorize derivation length, whereas 2. is applied in the very last step to obtain $\g'(0)$ as upper bound.

\smallskip

For what concerns Arai's Theorems 5 and 6, the situation is as follows. Theorem 6 allows to extract, from a derivation in $Diag(\emptyset)+(prg)^{\emptyset}$ of the well-foundedness of a linear relation $\prec$, an upper bound for the order type of $\prec$; since both Theorem 6 and its proof do not refer to \g, they remain untouched. Theorem 5, which plays a crucial role in Arai's proof, is a WP-elimination result that allows to remove the use of the rule $(WP)$ (which is a well-ordering principle in the form of a rule) in a formal derivation. More precisely, Theorem 5 states that, under some side conditions, from a derivation of a set of formulas $\Gamma$ in $(prg)^{\infty} + (WP)$ one can obtain a derivation of $\Gamma$ in $Diag(\emptyset)+(prg)^{\emptyset}$. 
The proof of Arai's Theorem 5 is preceded by two auxiliary lemmas, Lemma 1 and Lemma 2 \cite[pag. 272]{Arai20}, which again do not refer to \g\ and thus are untouched. Considering the proof of Theorem 5 itself, a key step is obtaining, from an \ita{embedding}, i.e., an injective order-preserving function, $f\! :\, <^{\emptyset}_A \frec \omega^{F(\ga,\al)+\ga+1}$, an embedding $\mathcal{F}\! :\, <^{\emptyset}_{\g(A)} \frec \g(\omega^{F(\ga,\al)+\ga+1})$ with $\g(\omega^{F(\ga,\al)+\ga+1})$ well-ordered by WOP(\g). %\footnote{By Theorem 6, from $\{n\}\vdash^{F(\ga, \al)+ \ga}E_{A}(n)$ for all $n$, we obtain an embedding $f\! :\, <^{\emptyset}_A \frec \omega^{F(\ga,\al)+\ga+1}$ and from this, with Proposition 1, an embedding $F\! :\, <^{\emptyset}_{\g(A)} \frec \g(\omega^{F(\ga,\al)+\ga+1})$.}. 
In Arai's proof, this step is ensured by \cite[Proposition 1]{Arai20} that relies on the term structure conditions imposed on $\g$, namely \cite[Definition 3]{Arai20}. Here instead, this fact is derived from the property of $\g$ of being weakly increasing: if $f$ is such an embedding, then $ot(<^{\emptyset}_{A}) \! < \! \omega^{F(\ga,\al)+\ga+1}$, thus $\g(ot(<^{\emptyset}_{A})) \! \leq \! \g(\omega^{F(\ga,\al)+\ga+1})$ and we can consider as $\mathcal{F}\! :\, <^{\emptyset}_{\g(A)} \frec \g(\omega^{F(\ga,\al)+\ga+1})$ the identity function\footnote{We warn the reader not to be baffled by the two distinct roles played here by the function letters $F$ and $\mathcal{F}$: binary ordinal function, $F(\al, \be)$, and embedding, $\mathcal{F}\! : \, <^{\emptyset}_{\g(A)} \frec \g(\omega^{F(\ga,\al)+\ga+1})$; we use the same letter to keep the notation as equal as possible to Arai's article.}. In conclusion, also Theorem 5 holds under our weaker hypotheses for $\g$ and so Arai's proof. \qed

\begin{rem}\label{rem:NotTerm} Differently from Arai, aside from being computable, we do not required for the term system for $\g$ any other property, such as, for example, the ones in \cite[Definition 3]{Arai20}; in particular, the final embedding $\mathcal{F}\! :\, <^{\emptyset}_{\g(A)} \frec \g(\omega^{F(\ga,\al)+\ga+1})$ may not be definable in \ACA. As a minor drawback, our proof is not directly formalizable in \ACA. Instead, a metatheory strong enough to treat directly ordinals, such as ZFC, is required.
\end{rem}

%\subsection{Ordinal analysis through a Well-Ordering rule}\label{subsec:WOPrule}

\subsection{Proof-Theoretic Ordinals of WOPs}\label{subsec:WOP}

The extension of Arai's theorem previously obtained is now applied to compute the proof-theoretic ordinals of WOP's related to specific ordinal functions. More precisely, we apply it to the first two of the following ordinal functions.

\begin{defi}\label{def:ordinalfun} Let us define: % the following ordinal functions:

\begin{enumerate}

\item $\g_{\omega}(\fX) \putAs \vartheta(\Omega^{\omega} \x \fX)$;

\item $\g_{\forall}(\fX) \putAs \sup_n \vartheta(\Omega^n \x \fX)$;

\item $\g_{\times}(\fX) \putAs \fX^{\omega}$; %\sup_n \fX^{n} =

\item $\g_{+}(\fX) \putAs  \fX \x \omega$. %\sup_n \fX \x n =

\end{enumerate}

\end{defi}
These ordinal functions rise, in this paper, from the analysis of some closure properties for \wqo: the former two are related to Kruskal's theorem, whereas the latter two to closure under products and disjoint unions; see respectively  \Subsec{subsec:KT} and \Subsec{subsec:wqo} for more details. Regarding the term structures needed to express these functions in second-order logic, the situation is as follows: for $\g_{\times}$ and $\g_{+}$, the corresponding term structures are relatively simple and can be found in the literature. In contrast, the cases of $\g_{\omega}$ and $\g_{\forall}$ are more complex; therefore, for the sake of brevity, we present the term structure only for $\g_{\omega}$. As many others, this term structure is definable in~\RCA.

Given a well-order $(X,\leq_X)$, we inductively build the well-order $(G_{\omega}(X), \leq_G~\!)$ defining simultaneously both the elements of $G_{\omega}(X)$ and the relation $\leq_G$:

\begin{enumerate}

\item $0 \in G_{\omega}(X)$, and, for each $x \in X \setminus \{0\}$, a constant $c_x$ belongs to $G_{\omega}(X)$;

\item if $x \in X$ and $\al_0, \dots, \al_n \in G_{\omega}(X)\setminus \{0\}$, then %$\vartheta(\Omega^{\omega} \x x) \in G_{\omega}(X)$ and 
$\vartheta(\Omega^{\omega} \x c_x + \Omega^n \x \al_n + \dots + \Omega^0 \x \al_0) \in G_{\omega}(X)$;

\item if $\al_0, \dots, \al_n \in G_{\omega}(X)$, $\al_i$ starts with $\vartheta$ for all $i \in \{0, \dots, n\}$, and $ \al_0 \geq_G \dots \geq_G \al_n$, then $\al_0 + \dots + \al_n \in G_{\omega}(X)$;

\item if $0 \neq \al$, then $0 <_G \al$, and $c_x \leq_G c_y$ if and only if $x \leq_X y$;

\item $\al_0 + \dots + \al_n <_G \be_0 + \dots +\be_m$ if and only if:
    \begin{enumerate}
    
      \item $n < m$ and $\forall i \! \leq \! n\, \al_i = \be_i$, or
      
      \item $\exists i \leq \min(n,m)$ such that $\forall j \! < \! i\, \al_j = \be_j \ \et \ \al_i <_G \be_i$;
    
    \end{enumerate}

\item $\al_0 + \dots + \al_n <_G \vartheta( \Omega^{\omega} \x c_x + \Omega^m \x \be_m + \dots + \Omega^0 \x \be_0) = \be$ if and only if $\forall i \! \leq \! n\ \al_i \! <_G \! \be$, with the side condition $n \geq 1$;

\item $\al = \vartheta( \Omega^{\omega} \x c_x + \Omega^n \x \al_n + \dots + \Omega^0 \x \al_0) <_G \be_0 + \dots + \be_m$ if and only if $\al <_G \be_0$, with the side condition $m \geq 1$;

\item $\al = \vartheta(\Omega^{\omega} \x c_x + \Omega^n \x \al_n + \dots + \Omega^0 \x \al_0) <_G \be = \vartheta(\Omega^{\omega} \x c_y + \Omega^m \x \be_m + \dots + \Omega^0 \x \be_0)$ if and only if:
     \begin{enumerate}
   
        \item $\al \leq_G \be_j$ for some $j \leq m$, or
        
        \item $x <_X y$ and $\forall i \! \leq \! n\ \al_i \! <_G \! \be$, or
        
        \item $x=y$, $n=m$, and $\exists i \! \leq \! n\, \forall j \! < \! i \, \al_j = \be_j \ \et \ \al_i \! <_G \! \be_i$, or
        
        \item $x=y$, $n <m$, and $\forall i \! \leq \! n \ \al_i \! <_G \! \be$.
   
     \end{enumerate}

\end{enumerate}
For what concerns ordinal analysis of these functions, we start with $\g_{\omega}$ and~$\g_{\forall}$. %We compute now the proof-theoretic ordinals of the WOPs corresponding to these four ordinal functions. More precisely, regarding  the WOPs of the first two, we calculate their proof-theoretic ordinals over ACA$_0$, thus resorting to \Theor{theor:WOP}; for the WOPs of the last two instead, since they are provable over ACA$_0$, we compute their proof-theoretic ordinals over RCA$_0$ using \Theor{theor:WOPrule}.

\begin{prop}\label{prop:ghyp} Let $\al, \be$ be countable ordinals. Both $\g_{\omega}(\fX)=\vartheta(\Omega^{\omega} \x \fX)$ and $\g_{\forall n}(\fX)=\sup_{n}\vartheta(\Omega^n \x \fX)$ satisfy the following properties:

\begin{enumerate}

\item $\al \! \leq \! \be$ implies $\g(\al) \! \leq \! \g(\be)$; % and 

\item $\omega^{\g'(0)}=\g'(0)$.

%\item for every strictly increasing function $f\! : \al \frec \be$ there exists a strictly increasing function $F\! : \g(\al) \frec \g(\be)$ such that $F(\g(\ga))=\g(f(\ga))$.

\end{enumerate}

\end{prop}

\noindent \ita{Proof:} Both properties derive directly from the definitions of \g$_{\omega}$ and \g$_{\forall n}$, together with the properties of the collapsing function $\vartheta$ (see \Lem{lem:vartheta}), e.g., $\vartheta(\al)$ is always an epsilon number, namely $\omega^{\vartheta(\al)}= \vartheta(\al)$. \qed % Property 5.2 is less evident and we prove it for the trickiest of the two cases, namely \g$_{\forall n}$. 

%\smallskip

%\tbe

\smallskip

Next, we compute $\g_{\omega}'(0)$ and $\g_{\forall n}'(0)$.

\begin{prop}\label{prop:g'0} $\g_{\omega}'(0)=\vartheta(\Omega^{\omega+1})$ and $\g_{\forall n}'(0)=\vartheta(\Omega^{\omega}+\omega)$.
\end{prop}

\noindent \ita{Proof:} We focus in particular on the second case, $\g_{\forall n}'(0)=\vartheta(\Omega^{\omega}+\omega)$, being the most involved. Let us define the following increasing sequence $\al_0=\omega$ and $\al_{k+1}=\sup_{n}\vartheta(\Omega^{n}\x \al_k)$; given the definition of $\vartheta$ (see \Subsec{subsec:preliminaries}), one can compute its first elements, i.e., $\al_1=\sup_{n}\vartheta(\Omega^n \x \omega)=\vartheta(\Omega^{\omega})$ and $\al_2=\sup_{n}\vartheta(\Omega^{n} \x \vartheta(\Omega^{\omega}))=\vartheta(\Omega^{\omega} + 1)$. We prove by induction on $k$, with the base case given by $\al_1$, that $\al_{k+1}=\vartheta(\Omega^{\omega}+k)$; from this $\g_{\forall n}'(0)=\sup_{k}\al_k=\vartheta(\Omega^{\omega}+\omega)$ easily follows. 

First we prove $\al_{k+1} \leq \vartheta(\Omega^{\omega}+k)$. So let $k\! \geq \! 1, \al_k=\vartheta(\Omega^{\omega}+k-1)$ and $\al_{k+1}=\sup_{n}\vartheta(\Omega^n \x \vartheta(\Omega^{\omega}+k-1))$; by definition of $\vartheta$, $\vartheta(\Omega^n \x \vartheta(\Omega^{\omega}+k-1)) \! < \! \vartheta(\Omega^{n+1} \x \vartheta(\Omega^{\omega}+k-1)) \! < \! \vartheta(\Omega^{\omega}+k)$ holds for all $n$ and thus $\al_{k+1}=\sup_{n}\vartheta(\Omega^n \x \vartheta(\Omega^{\omega}+k-1)) \! \leq \! \vartheta(\Omega^{\omega}+k)$. 
For the other direction, $\al_{k+1} \geq \vartheta(\Omega^{\omega}+k)$, we resort to an equivalent definition for $\vartheta$, namely $\vartheta(\al) \! = \! \min \{\xi \! \in \! E \, |\, {\al}^* \! < \! \xi  \ \mbox{and}\ \forall \be \! < \! \al \, ( \be^* \! < \! \xi \frec \vartheta(\be) \! < \! \xi) \}$ \cite[pag. 15]{VanderMereen15}, where $E$ is the set of epsilon numbers and, given an ordinal $\ga$, $\ga^*$ is the maximum epsilon number, apart from $\Omega$, in the Cantor normal form of $\ga$. Keeping in mind the aforementioned alternative definition for $\vartheta(\al)$, since $\al_{k+1}=\sup_{n}\vartheta(\Omega^n \x \vartheta(\Omega^{\omega}+k-1))$ is an epsilon number and $(\Omega^{\omega}+k)^*=0$, to prove $\al_{k+1} \geq \vartheta(\Omega^{\omega} + k)$ it remains to check that, for all $ \be \! < \! \Omega^{\omega}+k$, if $ \be^* \! < \! \al_{k+1}$, then $\vartheta(\be) \! < \! \al_{k+1}$. 
Let us assume that $\be \! < \! \Omega^{\omega}+k$ and $\be^* \! < \al_{k+1} = \! \sup_{n}\vartheta(\Omega^{n}\x \vartheta(\Omega^{\omega}+k-1))$, the first condition amounts to $\be \! < \! \Omega^{\omega}$ or $\be \! \in \! \{ \Omega^{\omega}, \Omega^{\omega}+1, \dots, \Omega^{\omega}+k-1\}$. The latter case case is the simplest one, if $\be \! \in \! \{ \Omega^{\omega}, \Omega^{\omega}+1, \dots, \Omega^{\omega}+k-1\}$, then $\be = \Omega^{\omega} + r$, with $0 \! \leq \! r \! \leq \! k-1$, and in this case $\vartheta(\be) \leq \vartheta (\Omega^{\omega}+k-1) = \al_k \leq \al_{k+1}$.
We consider now the former, if $\be \! < \! \Omega^{\omega}$ and $\be^* \! < a_{k+1}=\! \sup_{n}\vartheta(\Omega^n \x \vartheta(\Omega^{\omega}+k-1))$, then there exists a natural number $\bar{n}$ such that $\be \! < \! \Omega^{\bar{n}}\x \vartheta(\Omega^{\omega}+k-1)$ and $\be^* \! < \! \vartheta(\Omega^{\bar{n}} \x \vartheta(\Omega^{\omega}+k-1))$; but, since $\vartheta(\Omega^{\bar{n}} \x \vartheta(\Omega^{\omega}+k-1))$ is a $\vartheta$ number, then $\vartheta(\be) \! < \! \vartheta(\Omega^{\bar{n}} \x \vartheta(\Omega^{\omega}+k-1)) \! < \! \sup_{n}\vartheta(\Omega^n \x \vartheta(\Omega^{\omega}+k-1))=\al_{k+1}$. Thus $\al_{k+1}\geq \vartheta(\Omega^{\omega}+k)$ and so $\al_{k+1}=\vartheta(\Omega^{\omega}+k)$.

The case $\g_{\omega}'(0)=\vartheta(\Omega^{\omega+1})$ where $\g_{\omega}'(0) \putAs \min \{ \al \! > \! 0 \, | \, \forall \be \! < \! \al \ \vartheta( \Omega^{\omega} \x \be) \! < \! \al  \}$ is treated in an analogous way considering this time the sequence $\al_0=\omega$ and $\al_{k+1}=\vartheta(\Omega^{\omega}\x \al_k)$. \qed

\smallskip

The following computations are now immediate.

\begin{cor}\label{cor:WOPordinals} 
\[
|\mbox{\textsf{ACA}}_0 + \mbox{WOP}(\g_{\omega})|= \vartheta \! \left( \Omega^{\omega+1} \right)\ \,  \mbox{and} \ \ |\mbox{\textsf{ACA}}_0 + \mbox{WOP}(\g_{\forall})|= \vartheta \! \left( \Omega^{\omega} + \omega \right).
\]
\end{cor}

\proof We apply together \Theor{theor:WOP}, \Prop{prop:ghyp} and \Prop{prop:g'0}.~$\square$ %\qed

\medskip

We move to the other two ordinal functions, namely $\g_{\times}(\fX)=\fX^{\omega} $ and $\g_{+}(\fX)=\fX \x \omega$. Since WOP$(\g_{\times})$ and WOP$(\g_{+})$ are provable over \ACA, with WOP$(\g_{+})$ already provable over \RCA\, \cite{Hirst94,Hirst98}, we consider as base two weaker theories: \RCA\, for WOP$(\g_{\times})$ and \RCAs\, for WOP$(\g_{+})$. %As the reader may have noticed, \Theor{theor:WOP} uses as base theory ACA$_0$. Although an equivalent result for RCA$_0$ seems achievable, such a general computation is not already available; thus, in order to obtain at least a preliminary estimation relatively to these ordinal function with a Well Ordering Principle (WOP) but with a \ita{Well Ordering Rule}(WOR).

First of all, we compare WOP$(\g_{\times})$ and WOP$(\g_{+})$ with two uniform WOP's strictly correlated; more precisely, let us consider the following statements:
\begin{itemize}

	\item $\forall n\, \mbox{WOP}(\times^n):\ \forall n\, \forall \fX\, [\mbox{WO}(\fX) \frec \mbox{WO}(\fX^n)]$;
	
	\item $\forall n\, \mbox{WOP}(+^n):\ \forall n\, \forall \fX\, [\mbox{WO}(\fX) \frec \mbox{WO}(\fX \x n)]$.
	
\end{itemize}

\begin{prop}\label{prop:equiv} The following equivalences hold:

\begin{enumerate}

  \item \textsf{RCA}$^*_0 \, \vdash \ \mbox{WOP}(g_{\times}) \, \leftrightarrow \, \forall n\, \mbox{WOP}(\times^n)$;
  
  \item \textsf{RCA}$^*_0 \, \vdash \ \mbox{WOP}(g_{+})\, \leftrightarrow \, \forall n\, \mbox{WOP}(+^n)$.
  
\end{enumerate}

\end{prop}
\proof We consider only the first one which derives from $\sup_n(\al^n)=\al^\omega$ and the following chain of equivalences
\[
 \forall n\, \mbox{WOP}(\times^n) \, \leftrightarrow \, \forall \fX\, [\mbox{WO}(\fX) \frec \forall n\,\mbox{WO}(\fX^n)] \, \leftrightarrow \, \forall \fX\, [\mbox{WO}(\fX) \frec \mbox{WO}(\sup_n(\fX^n))].
\]
The second point is treated analogously. \qed

\noindent The difference with respect to WOP$(\g_{\omega})$ and WOP$(\g_{\forall})$ is due to the fact that $\sup_n(\al \x n)=\al \x \omega$ and $\sup_n(\al^n)=\al^{\omega}$, whereas $\sup_n (\vartheta(\Omega^n \x \al)) \neq  \vartheta(\Omega^{\omega} \x \al)$.\footnote{A simple counterexample is given by $\al= \Omega$.}

\smallskip

We now analyse WOP$(\g_{\times})$ and WOP$(\g_{+})$. For WOP$(\g_{\times})$, there are the following two results, where $\Sigma^0_2\mbox{-IND}$ denotes the induction schema:
\[
\fhi(0) \et \forall n\, (\fhi(n) \frec \fhi(n+1)) \ \frec \ \forall n\, \fhi(n),
\]
for all $\Sigma^0_2$ formulas.

\begin{lem}\label{lem:Patrick} \RCA $\, \vdash \ \mbox{WOP}(\g_{\times}) \ \longleftrightarrow \ \Sigma^0_2\mbox{-IND}$.
\end{lem}
\proof This has been proved by Patrick Uftring \cite{Uftring23}(see also \cite{FU23}). \qed

\begin{prop}\label{prop:|gX|} $|\mbox{\textsf{RCA}}_0 + \mbox{WOP}(\g_{\times})|=\omega^{\omega^{\omega}}$.
\end{prop}

\proof Directly from \Lem{lem:Patrick} and the ordinal analysis of $\Sigma^0_2\mbox{-IND}$. \qed

\smallskip 

Regarding WOP$(\g_{+})$ instead, the next proposition holds.

\begin{prop}\label{prop:|g+|} $|\mbox{\textsf{RCA}}_0^* + \mbox{WOP}(\g_{+})|=\omega^{\omega}$.

\end{prop}

\proof The proof is given by the following series of inequalities:
\[
\omega^{\omega} \leq |\mbox{\textsf{RCA}}_0^* + \mbox{WOP}(\g_{+})| \leq |\mbox{\textsf{RCA}}_0| \leq \omega^{\omega}.
\]
For the first inequality, we apply the same strategy used at the beginning of the proof of \Theor{theor:WOP}, relying this time on the sequence $\al_0 \putAs \omega$ and $\al_{k+1}\putAs \al_k \x \omega$ obtaining that $|\mbox{\textsf{RCA}}_0^* + \mbox{WOP}(\g_{+})| \geq \sup_k \al_k = \omega^{\omega}$. The second inequality derives from the fact that both \RCAs\, and WOP$(\g_{+})$ are provable in \RCA, and the last one from the ordinal analysis of \RCA. \qed

%%%%%%%
%%%%%%%

\section{Well Quasi-Orders Closure Properties}\label{sec:wqo}

In this section, previous computations are applied to calculate the proof-theoretic ordinals of some {\wqo} closure properties, namely results of the form: \vir $Q$ {\wqo} implies $\sF(Q)$ {\wqo}'' with $\sF$ some set operation; the similarity with \eqref{eq:WOP} should be apparent. We consider two applications: the former regards two different versions of labelled Kruskal's theorem (see \Theor{theor:Kruskal}), with the unlabelled version already treated in \cite{RW93}; the latter concerns some closure properties related to disjoint unions and products.

\subsection{Proof-theoretic Ordinals of Kruskal's Theorem}\label{subsec:KT}

We compute now, using \Theor{theor:WOP}, the aforementioned estimations for the proof-ordinals of Kruskal's theorem with labels, namely:
\[
|\mbox{\textsf{RCA}}_0 + \forall n\, \mbox{KT}_{\ell}(n)|=\vartheta (\Omega^{\omega}+\omega)\ \mbox{and}\ |\mbox{\textsf{RCA}}_0 + \mbox{KT}_{\ell}(\omega)|=\vartheta (\Omega^{\omega+1}).
\] 
Instrumental to our goal, we prove the following equivalences:

\begin{theor}\label{theor:KTl-WOP} {(\RCA)} The following are equivalent:

\begin{enumerate}

\item KT$_\ell (\omega) \! : \, \forall Q\, [Q\, \mbox{{\wqo}} \frec \mT (Q)\, {\wqo}]$;

\item WOP$(\g_{\omega})\! : \, \forall \fX\, [ \mbox{WO}(\fX) \frec \mbox{WO}(\vartheta(\Omega^{\omega}\x \fX)) ]$. 

\end{enumerate}
For each $n$, the following are equivalent:

\begin{enumerate}

\item KT$_\ell (n) \! : \, \forall Q\, [Q\, \mbox{{\wqo}} \frec \mT_n (Q)\, {\wqo}]$;

\item $\forall \fX\, [ \mbox{WO}(\fX) \frec \mbox{WO}(\vartheta(\Omega^{n}\x \fX)) ]$. 

\end{enumerate}

\end{theor}

\noindent \ita{Proof}: We consider the equivalence between the first two points reasoning in \ACA, since they both imply \ACA\, over \RCA.

$1) \imp 2)$ Let $\fX$ be a well-order, then $\fX$ is a {\wqo} and, by $1)$, also $(\mT(\fX), \preceq)$ is a {\wqo}; in particular any extension of the tree embeddability $\preceq$ is well-founded. But, since $\preceq $ can be extended over $\mT(\fX)$ to a linear order order-isomorphic to $\vartheta(\Omega^{\omega}\x \fX)$ \cite{RW93}, $\vartheta(\Omega^{\omega}\x \fX)$ is a well-founded linear order and thus is well-ordered. 

 $2) \imp 1)$ Let $Q$ be a {\wqo}, then the tree $T_B$ of bad sequences in $Q$ is well-founded and, provable in \ACA\, \cite[lemma V.1.3]{Simpson09}, the Kleene-Brouwer linearization $\fX$ of $T_B$ is a well-order. Thus $Q$ admits a \ita{reification} by the well-order $\fX$, namely there is a function $f \! : \, T_B \frec \fX +1$ such that if the bad sequence $a$ extends the bad sequence $b$, i.e., $b$ is an initial segment of $a$, then $f(a) \leq f(b)$. By \cite[Theorem 2.2]{RW93}, this implies the existence of a reification of $\mT(Q)$ by $\vartheta(\Omega^{\omega} \x \fX)$ which, thanks to $2)$, is a well-order; thus $\mT(Q)$ is a {\wqo} since it has no infinite bad sequences.
 
\noindent For each $n$, the second part is proved analogously. \qed

\smallskip
\noindent From the second equivalence of the previous theorem, one can easily obtain:

\begin{prop}\label{prop:KTln-WOP}{(RCA$_0$)} The following are equivalent:

\begin{itemize}

\item $\forall n\,$KT$_\ell (n)$ \!: $\forall n\, \forall Q\, [Q\, \mbox{{\wqo}} \frec \mT_n (Q)\, {\wqo}]$;

\item $\forall n\, \forall \fX\, [ \mbox{WO}(\fX) \frec \mbox{WO}(\vartheta(\Omega^{n}\x \fX)) ]$;

\item $\forall \fX\, [ \mbox{WO}(\fX) \frec \forall n\, \mbox{WO}(\vartheta(\Omega^{n}\x \fX)) ]$;

\item $\forall \fX\, [ \mbox{WO}(\fX) \frec \mbox{WO}(\sup_n (\vartheta(\Omega^{n}\x \fX))) ]$.

\end{itemize}

\end{prop}
\noindent \ita{Proof:} Straightforward using \Theor{theor:KTl-WOP} (cf. \Prop{prop:equiv}). \qed

Finally, our main result is derived:

\begin{theor}\label{theor:KTlordinals}
\[
|\mbox{\textsf{RCA}}_0 + \forall n\, \mbox{KT}_{\ell}(n)|=\vartheta (\Omega^{\omega}+\omega)\ \mbox{and}\ |\mbox{\textsf{RCA}}_0 + \mbox{KT}_{\ell}(\omega)|=\vartheta (\Omega^{\omega+1}).
\] 
\end{theor}

\noindent \ita{Proof:} We simply apply together \Cor{cor:WOPordinals}, \Theor{theor:KTl-WOP}, \Prop{prop:KTln-WOP}.

\smallskip

For sake of completeness, the analogous result for the unlabelled case is reported; followed by some corollaries.

\begin{theor}\label{theor:KT} \RCA$\, \vdash \ \forall n\, \mbox{KT}(n) \sse \mbox{KT}(\omega) \sse \mbox{WO}(\vartheta(\Omega^{\omega}))$.
\end{theor}

\noindent \ita{Proof:} \cite{RW93}. \qed

\begin{cor}\label{theor:KTln-/-KTlw} \RCA$ \ \nvdash \ \forall n\, \mbox{KT}_{\ell}(n) \ \frec \ \mbox{KT}_{\ell}(\omega)$.
\end{cor}

\noindent \ita{Proof:} Directly from \Theor{theor:KTlordinals}. \qed

\begin{cor}$|\mbox{\textsf{RCA}}_0 + \mbox{KT}(\omega)|=\vartheta(\Omega^{\omega})^{\omega}$.
\end{cor}

\noindent \ita{Proof:} Combining together \Theor{theor:KT} and \cite[Theorem 5]{CMR19}. \qed

\smallskip

We conclude this section with a general equivalence result between closure under well quasi-orders and well-ordering principles. Let $\sG$ be a set operation that preserves {\wqo}, e.g. the string operator $(Q \mapsto Q^*)$ of Higman's lemma \cite{Higman52}, and, for an ordinal $\al$, let us define the following ordinal function $\g(\al) \putAs ot(\sG(\al))$, i.e., the order type of $\sG(\al)$ which is a {\wqo} given the preservation property of $\sG$; then the following holds.

\begin{theor}\label{theor:Equiv} Let $\sG$ and $\g$ be as before. Assume that, for every quasi-order $Q$, \ACA\, suffices to prove that 
\begin{itemize}

\item if $Q$ has order type $\al$, then $\sG(Q)$ has order type $\g(\al)$,

\item if $Q$ has a reification by $\al$, then $\sG(Q)$ has a reification by $\g(\al)$.

\end{itemize}
Then the following are equivalent over \ACA:

\begin{enumerate}

\item $\forall Q\, [Q\, \mbox{{\wqo}} \frec \sG (Q)\, {\wqo}]$;

\item $\forall \fX\, [ \mbox{WO}(\fX) \frec \mbox{WO}(\g(\fX)) ]$. 

\end{enumerate}

\end{theor}

\proof The idea is to generalize the proof of \Theor{theor:KTl-WOP} working over \ACA.

$1) \imp 2)$ Let $\fX$ be a well-order, then $\fX$ is a {\wqo} of order type $\fX$ and, by $1)$ and the order type hypothesis for \ACA, $\sG(\fX)$ is a {\wqo} of order type $\g(\fX)$. Since $\sG(\fX)$ is {\wqo}, this means that $\g(\fX)$ is a well-founded linear extension of $\sG(\fX)$ and thus is a well-order.

$2) \imp 1)$ Let $Q$ be a {\wqo}, then the tree $T_B$ of bad sequences in $Q$ is well-founded and, provable in \ACA, the Kleene-Brouwer linearization $\al$ of $T_B$ is a well-order. Thus $Q$ has a reification by $\al$ and $\sG(Q)$ has a reification by $\g(\al)$ which, thanks to $2)$, is a well-order; thus $\sG(Q)$ is a {\wqo}. \qed

The previous theorem can be easily extended to the case of a family of set operations with a natural parameter, such as in \Prop{prop:KTln-WOP}.

\subsection{Wqo Closure Properties over \RCA\, and \RCAs}\label{subsec:wqo}

We investigate now three other closure properties for well quasi-orders related to product, sum and disjoint union of {\wqo}'s. First, let us recall how these quasi-orders constructions are defined.
\begin{defi}\label{def:prod}
Given two qo $(P, \leq_P)$ and $(Q, \leq_Q)$, we define:
\begin{enumerate}

\item the product $(P \times Q, \leq_{\times})$ with 
\[
(p_1, q_1) \leq_{\times} (p_2, q_2) :\Longleftrightarrow p_1 \leq_P p_2 \, \et \, q_1 \leq_Q q_2;
\]

\item the disjoint union $(P \, \dot{\cup} \, Q, \leq_{\dot{\cup} })$ with
\[
p \leq_{\dot{\cup}} q :\Longleftrightarrow  (p,q \in P \et p \leq_P q) \, \vel \, (p,q \in Q \et p \leq_Q q);
\]

\item the sum\footnote{Sum and disjoint union of qo's have the same base set, $P \, \dot{\cup} \, Q$, they differ only for the quasi-order relation; in particular the sum qo $\leq_{+}$ is an extension of the disjoint qo $\leq_{\dot{\cup}} $.}
$(P \, \dot{\cup} \, Q, \leq_{+ })$, denoted by $P+Q$, with
\[
p \leq_{+} q :\Longleftrightarrow (p \! \in \! P \et q \! \in \! Q) \, \vel \, p \leq_{\dot{\cup}} q.
\]

\end{enumerate}
\end{defi}

Next, we present the closure properties under consideration.
\begin{defi}For what concerns product of {\wqo}, the following principles are considered:
\begin{enumerate}

 \item WQP$(\times):\ \forall P, Q \ [P, Q \, \mbox{{\wqo}} \, \frec \, P \times Q \, \mbox{{\wqo}}]$;
 
 \item WQP$(\times^n):\ \forall Q\ [Q\, \mbox{{\wqo}} \, \frec \, \underbrace{Q \times \dots \times Q}_{n \, \mbox{times}}\, \mbox{{\wqo}}]$;
 
 \item WQP$(\times^{\omega}) \putAs \forall n\, \mbox{WQP}(\times^n):\ \forall n\, \forall Q\ [Q\, \mbox{{\wqo}} \, \frec \, \underbrace{Q \times \dots \times Q}_{n \, \mbox{times}}\, \mbox{{\wqo}}]$.
\end{enumerate}

\noindent For what concerns sum of {\wqo}, the following principles are considered:
\begin{enumerate}

 \item WQP$(+):\ \forall P, Q \ [P, Q \, \mbox{{\wqo}} \, \frec \, P + Q \, \mbox{{\wqo}}]$;
 
 \item WQP$(+^n):\ \forall Q\ [Q\, \mbox{{\wqo}} \, \frec \, \underbrace{Q + \dots + Q}_{n \, \mbox{times}}\, \mbox{{\wqo}}]$;
 
 \item WQP$(+^{\omega}) \putAs \forall n\, \mbox{WQP}(+^n):\ \forall n\, \forall Q\ [Q\, \mbox{{\wqo}} \, \frec \, \underbrace{Q+ \dots + Q}_{n \, \mbox{times}}\, \mbox{{\wqo}}]$.
\end{enumerate}

\noindent For what concerns disjoint union of {\wqo}, the following principles are considered:
\begin{enumerate}

 \item WQP$(\dot{\cup}):\ \forall P, Q \ [P, Q \, \mbox{{\wqo}} \, \frec \, P \, \dot{\cup} \, Q \, \mbox{{\wqo}}]$;
 
 \item WQP$(\dot{\cup}^n):\ \forall Q\ [Q\, \mbox{{\wqo}} \, \frec \, \underbrace{Q \, \dot{\cup} \, \dots \, \dot{\cup} \, Q}_{n \, \mbox{times}}\, \mbox{{\wqo}}]$;
 
 \item WQP$(\dot{\cup}^{\omega}) \putAs \forall n\, \mbox{WQP}(\dot{\cup}^n):\ \forall n\, \forall Q\ [Q\, \mbox{{\wqo}} \, \frec \, \underbrace{Q \, \dot{\cup} \, \dots \, \dot{\cup} \, Q}_{n \, \mbox{times}}\, \mbox{{\wqo}}]$.
\end{enumerate}

\end{defi}
For sake of readability, from now on we generally omit the under brackets with the number of iterations which will be clear from the context. Before treating the ordinal analysis of WQP$(\times^{\omega})$, WQP$(+^{\omega})$, and WQP$(\dot{\cup}^{\omega})$, we recall and state some general results regarding the provability, and the proof-theoretic relations, of WQP$(\times)$, WQP$(+)$, and WQP$(\dot{\cup})$.

\begin{prop}\label{prop:wqp} The following hold:

\begin{enumerate}

 \item \textsf{RCA}$_0 \, \nvdash \ \mbox{WQP}(\times) $;
 
 \item \textsf{RCA}$^*_0 \, \vdash \ \mbox{WQP}(+) $;
 
 \item \textsf{RCA}$^*_0 \, \vdash \ \mbox{WQP}(\dot{\cup}) $.

\end{enumerate}

\end{prop}

\proof \ita{1.} can be actually strengthened to \textsf{WKL}$_0 \, \nvdash \ \mbox{WQP}(\times) $, see \cite[Corollary 4.7]{CMS04}. 

\noindent \ita{2.} derives immediately from \ita{3.} since $p \leq_{\dot{\cup}} q$ implies $p \leq_{+} q$, and thus any bad sequence in $\leq_{+}$ gives rise to a bad sequence in $\leq_{\dot{\cup}}$.

\noindent \ita{3.} let us consider the finite pigeonhole principle RT$^1_2$: ``$\forall f \!: \mN \frec \{0,1\}$ there is $n \! \in \! \{0,1\}$ such that $f^{-1}(n)$ is infinite''. Hirst proved \cite[Theorem 6.3]{Hirst87} \textsf{RCA}$_0 \vdash \, \mbox{RT}^1_2$, and his proof suffices to obtain also \textsf{RCA}$^*_0 \vdash \, \mbox{RT}^1_2$. If $P,Q$ are {\wqo}, by applying RT$^1_2$ to an infinite bad sequence in $Q \, \dot{\cup} \, P$, we can extract an infinite bad subsequence in $P$ or $Q$, reaching a contradiction. \qed

\begin{prop}\label{prop:Xwqp} The following hold:

\begin{enumerate}

 \item \textsf{RCA}$^*_0 \, \vdash \ \mbox{WQP}(\times) \, \leftrightarrow \, \mbox{WQP}(\times^2)$;
 
 \item \textsf{RCA}$^*_0 \, \vdash \ \mbox{CAC} \, \frec \, \mbox{WQP}(\times) \, \frec \, \mbox{ADS} $.

\end{enumerate}

\end{prop}

\proof \ita{1.} ``$\frec$'' is immediate. The other direction derives from the fact that \textsf{RCA}$^*_0  \vdash \, \mbox{WQP}(\dot{\cup})$ and $P \times Q \subseteq (P \, \dot{\cup} \, Q) \times (P \, \dot{\cup} \, Q)$. 

\noindent \ita{2.} See \cite[Theorem 1.16]{Towsner20}, in particular both implications are strict. \qed

Applying together \Prop{prop:Xwqp} and \Cor{cor:CAC}, we immediately obtain the following.

\begin{cor} $|\mbox{RCA}_0 + \mbox{WQP}(\times)|=\omega^{\omega}$.

\end{cor}

For what concerns the ordinal analysis of WQP$(\times^\omega)$, WQP$(+^\omega)$, and WQP$(\dot{\cup}^\omega)$, since all three are provable over \ACA\, (see  for example \Theor{theor:ACAHigman}), in this last paragraph we consider \RCA, or sometimes \RCAs, as base theory.

\subsubsection{Closure under iterated product}

In order to obtain an ordinal analysis of WQP$(\times^\omega)$, a possible intermediate goal is to find a WOP$(*)$ equivalent to WQP$(\times^\omega)$ and then apply previous computations to estimate $|\mbox{\textsf{RCA}}_0 + \mbox{WOP}(*) |$. In this search we can take advantage from the calculations of de Jongh and Parikh \cite{dJP77} regarding the order type of the results of many operations on {\wqo}; in particular we have the following, where $\otimes$ is the natural product between ordinals \cite{Bachmann55}:

\begin{prop}\label{prop:prod} Given two {\wqo} $P$ and $Q$ with order type respectively $ot(P)=\al$ and $ot(Q)=\be$, then $ot(P\times Q)=\al \otimes \be$.
\end{prop}

Using this result as guide, if $Q$ as order type $\al$, then $Q \times ... \times Q$ has order type $\al^n$ (given $\al \otimes \al = \al \times \al$), and since we are consider this property uniformly on $n$ (``$\forall n\, \forall Q [\dots]$''), our candidate as ordinal function is $\g_{\times}(\al)=\sup_n \al^n=\al^{\omega}$; namely the WOP candidate is $\forall \fX \, [\mbox{WO}(\fX) \frec \mbox{WO}(\fX^{\omega})]$ already treated in \Subsec{subsec:WOP} (cf. \Prop{prop:equiv}).

The last step would be to prove the equivalence between WQP$(\times^\omega)$ and WOP$(\g_{\times})$ over \RCA; however, the next negative result holds.

\begin{prop}\label{prop:prodrel} Over \RCA, the following is the only provable implication between WQP$(\times^\omega)$ and WOP$(\g_{\times})$:
\[
\mbox{WQP}(\times^\omega) \ \longrightarrow \ \mbox{WOP}(\g_{\times}).
\]
\end{prop} 

\proof We consider separately each direction.

\smallskip

\noindent ``$\mbox{WQP}(\times^\omega)  \rightarrow  \mbox{WOP}(\g_{\times})$'' Let $\al$ be a well-order, we have to prove that $\al^{\omega}$ is a well-order; for proving this over \RCA, it suffices to prove that $\al^n$ is a well-order for every $n$. So let $n$ be a natural number, $\al$ is a well-order and thus a {\wqo} and, by WQP$(\times^\omega)$, also $\al \times \dots \times \al$ is a {\wqo}. But the standard lexicographic linear order $\preceq$ over $\al^n$ is a linear extension of the product order $\leq_{\times}$ over $\al \times \dots \times \al$ and, provably over \RCA, this implies its well-foundedness; thus $\al^n$ is a well-order.

%\noindent ``$\mbox{WOP}(\g_{\times})  \leftrightarrow \Sigma^0_2\mbox{-IND}$'' see \Lem{lem:Patrick}.%This has been proved by Patrick Uftring \cite{Uftring23}(see also \cite{FU23}).

\smallskip

\noindent ``$\mbox{WOP}(\g_{\times}) \ \slash \! \! \! \! \! \!  \rightarrow  \mbox{WQP}(\times^\omega) $'' Assume, by contradiction, that WOP$(\g_{\times})$ implies WQP$(\times^\omega)$ over \RCA; then, by \Lem{lem:Patrick}, also $\Sigma^0_2\mbox{-IND}$ implies WQP$(\times^\omega)$ which in turn, by \Prop{prop:Xwqp}, implies the ADS principle. But now we have reached a contradiction since in REC, the minimal $\omega$-model of \RCA, $\Sigma^0_2\mbox{-IND}$ holds but ADS does not. \qed

Given this negative result, to obtain the proof-theoretic ordinal of WQP$(\times^\omega)$, we rely on the following lemma (originally proved by Lorenzo Carlucci, see \cite{CMZ25}), and the subsequent characterization of RT$^2_{<\infty}$ due to Slaman and Yokoyama \cite{SY:pairs}.

\begin{lem}\label{lem:RT2} \RCA $\, + \mbox{RT}^{\, 2}_{<\infty} \, \vdash \, \mbox{WQP}(\times^\omega)$. 
\end{lem}

\proof Let us assume that there exist $n \! \in \! \mN$ and a {\wqo} $Q$ such that  $Q^n=\underbrace{Q \times \dots \times Q}_{n \, \mbox{times}}$ is not a {\wqo}.  Let $((q^i_1, \dots, q^i_n))_{i \in \mN}$ be a bad sequence in $Q^n$; then, for all $i<j$ we have $(q^i_1, \dots, q^i_n) \not \leq_{\times} (q^j_1, \dots, q^j_n)$, i.e., there exists $k \! \in \! \{0, \dots, n-1\}$ such that $q^j_k \not \leq_Q q^i_k$.

Define $c\! : [\mN]^2 \frec \{0, \dots, n-1\}$ as:
\[
c(i,j)\putAs \min\{k \! \in \! \{0, \dots, n-1\} \, | \, q^i_k \not \leq_Q q^j_k \}.
\]
Let $H$ be an infinite $c$-homogeneous set as given by RT$^2_{<\infty}$, and let $k^* \! \in \! \{0, \dots, n-1\}$ be the colour of $[H]^2$ under $c$. It is easy to see that $(q^h_{k^*})_{h \in H}$ is a bad sequence in $Q$, contradicting our hypothesis. \qed

\begin{theor}\label{theor:RT2} \textsf{WKL}$_0 + \mbox{RT}^{\, 2}_{<\infty}$ is $\Pi^1_1$-conservative over \RCA $\, + \mbox{B}\Sigma^0_3$.
\end{theor}

\proof See \cite[Theorem 2.1]{SY:pairs}. \qed

Combining the previous results with the ordinal analysis of $B\Sigma^0_n$ conducted in \Sec{sec:Bsigma}, we obtain the following.\footnote{The authors thank Lorenzo Carlucci for proposing this strategy.}

\begin{prop} $|\mbox{\textsf{RCA}}_0 + \mbox{WQP}(\times^\omega)|=\omega^{\omega^{\omega}}$.
\end{prop}

\proof The proof is given by the following series of inequalities:
\[
\omega^{\omega^{\omega}} \leq |\mbox{\textsf{RCA}}_0 + \mbox{WQP}(\times^\omega)| \leq |\mbox{\textsf{RCA}}_0 + \mbox{RT}^2_{<\infty}| = |\mbox{\textsf{RCA}}_0 + \mbox{B}\Sigma^0_3| = \omega^{\omega^{\omega}},
\]
where the first inequality comes from  \Prop{prop:|gX|} and \Prop{prop:prodrel}, and the second inequality from \Lem{lem:RT2}; while the first equality derives from \Theor{theor:RT2} and the second one from the case $n=3$ of \Cor{cor:Bsigma}.~$\square$% \qed

\subsubsection{Closure under iterated sum or disjoint union}

We consider now the principles WQP$(+^\omega)$ and WQP$(\dot{\cup}^\omega)$. The definitions of sum and disjoint quasi-orders, \Def{def:prod}, together with the fact that WQP$(+^\omega)$ and WQP$(\dot{\cup}^\omega)$ involve an arbitrary, but finite, number of copies of the same qo, allow to restate the two well quasi-orders principles under scrutiny in the following way, where we %identify a natural number with the set of its predecessors\footnote{Namely, $n \putAs \{ 0, 1, \dots, n-1\}$.} and 
add a third (but strictly related) principle.

\begin{defi}\label{def:altWQP} WQP$(+^\omega)$ and WQP$(\dot{\cup}^\omega)$ can be restate as:

\begin{itemize}

\item WQP$(+^\omega): \ \, \forall n\, \forall Q\, [Q \ \mbox{\wqo} \frec (n \times Q,\leq_{+})\ \mbox{\wqo}] $;

\item WQP$(\leq_{\times}^\omega)\!: \ \, \forall n\, \forall Q\, [Q \ \mbox{\wqo} \frec (n \times Q,\leq_{\times})\ \mbox{\wqo}] $;

\item WQP$(\dot{\cup}^\omega): \ \, \forall n\, \forall Q\, [Q \ \mbox{\wqo} \frec (n \times Q,\leq_{\dot{\cup}})\ \mbox{\wqo}] $.

\end{itemize}

\end{defi}
For sake of clarity, and to ease comparisons, we spell out the three aforementioned relations\footnote{Which are respectively, the lexicographic order $\leq_{+}$, the product order $\leq_{\times}$ and the disjoint order $\leq_{\dot{\cup}}$.} for given $m_1, m_2 < n$ and $q_1,q_2 \in Q$:

\begin{itemize}

\item $(m_1, q_1) \leq_{+} (m_2, q_2) :\Leftrightarrow m_1 < m_2 \, \vel \, (m_1 = m_2 \et q_1 \leq_Q q_2)$; 

\item $(m_1, q_1) \leq_{\times} (m_2, q_2) :\Leftrightarrow m_1 \leq m_2 \, \et \, q_1 \leq_Q q_2$;

\item $(m_1, q_1) \leq_{\dot{\cup}} (m_2, q_2) :\Leftrightarrow m_1 = m_2 \, \et \, q_1 \leq_Q q_2$.

\end{itemize}

We start the ordinal analysis of these principles with the next result.

\begin{lem}\label{lem:mXQ} The following hold:

\begin{enumerate}

\item \RCAs $ \, \vdash \ \mbox{WQP}(\dot{\cup}^\omega) \ \frec \ \mbox{WQP}(\leq_{\times}^\omega) \ \frec \mbox{WQP}(+^\omega)$;

\item \RCA $ \, \vdash \ \mbox{WQP}(+^\omega)$.

\end{enumerate}

\end{lem}

\proof \ita{1.}  $\leq_{+}$ is an extension, a linear extension actually, of $\leq_{\times}$ which, in turn, is an extension of $\leq_{\dot{\cup}}$; thus, a bad sequence in $(n \times Q,\leq_{+})$ would give rise to a bad sequence in $(n \times Q,\leq_{\times})$. Similarly between $(n \times Q,\leq_{\times})$ and $(n \times Q,\leq_{\dot{\cup}})$.

\smallskip

\noindent \ita{2.} Reasoning by contradiction, let us assume that there exist $m$ and $Q$ such that $Q$ is {\wqo}, but $(m \times Q,\leq_{+})$ is not. Given a bad sequence $((k_i, q_i))_i$ in $(m \times Q,\leq_{+})$, $(k_i)_i$ is a weakly descending sequence in $\{0, \dots, m-1\}$; so let $k$ be the minimum of $(k_i)_i$ and $r$ the first index for which $k_r=k$. Since $((k_i,q_i))_{i \geq r}$ is a bad sequence with $k_i=k$ for all $i \geq r$, $(q_i)_i$ must be a bad sequence in $Q$; which is impossible by hypothesis. \qed

\smallskip

WQP$(\dot{\cup}^\omega)$, similarly to WQP$(\leq_{\times}^\omega)$ but differently from WQP$(+^\omega)$, is not provable in \RCA; the unprovability derives directly from the next result together with the independence of the Infinite Pigeonhole Principle RT$^1_{<\infty}$ from RCA$_0$.

\begin{theor}\label{theor:4equiv} Over \RCA\, the following principles are equivalent:

\begin{enumerate}

\item WQP$(\dot{\cup}^\omega): \ \, \forall n\, \forall Q\, [Q \ \mbox{\wqo} \frec (n \times Q,\leq_{\dot{\cup}})\ \mbox{\wqo}] $;

\item WQP$(\leq_{\times}^\omega)\!: \, \forall n\, \forall Q\, [Q \ \mbox{\wqo} \frec (n \times Q,\leq_{\times})\ \mbox{\wqo}] $;

\item RT$^{\, 1}_{<\infty} :\ \forall n\, \forall f\! : \mN \frec \{0, \dots, n-1\} \ \exists A \in [\mN]^{\omega} \, \exists i < n\ \forall j \! \in \! A \ (f(j)=i)$;

\item B$\Sigma^0_2 : \forall x \leq a\, \exists y\, \fhi(x,y) \, \frec \, \exists z\, \forall x \leq a\, \exists y \leq z\, \fhi(x,y)$ for all $\fhi \in \Sigma^0_2$.

\end{enumerate}

\end{theor}

\proof The equivalence between RT$^{\, 1}_{<\infty}$ and B$\Sigma^0_2$ is well known in the literature (e.g. \cite{CJS01}); thus, we consider only the other three properties proving separately \ita{1.} $\imp$ \ita{2.}, \ita{2.} $ \imp$ \ita{3.} and \ita{3.} $ \imp $ \ita{1.} 

\smallskip

\noindent ``\ita{1.} $\imp$ \ita{2.}'' point \ita{1.} of \Lem{lem:mXQ}.

\smallskip

\noindent ``\ita{2.} $\imp$ \ita{3.}'' Assume by contradiction that RT$^1_{<\infty}$ fails; thus there exist a sequence $(n_i)_i$ with $i \! \in \! \mN$ and a bound $m \! \in \! \mN$ such $n_i < m$ for all $i \! \in \! \mN$, but no number less than $m$ occurs infinite times in $(n_i)_i$. We will construct a well-order $\al$ such that $(m \times \al, \leq_{\times})$ is not {\wqo}, contradicting WQP$(\leq_{\times}^\omega)$.\footnote{We adapt here the first part of the proof of \cite[Lemma 5.3.4.]{UftringPhD}.}

\noindent Let $\al$ be defined on $\mN$ as follows: For any two elements $i, j \in \mN$, we have $i <_{\al} j$ if and only if $n_i < n_j$ or both $n_i = n_j$ and $i >_{\mN} j$ hold. It can easily be seen that $\al$ is a linear order. For showing that $\al$ is also well-founded, consider a hypothetical infinite descending sequence $(i_k)_k$ in $\al$. By $\Pi^0_1$-induction, we know that there is a smallest $n < m$ such that there exists an index $\bar{k} \! \in \! \mN$ with $n = n_{i_{\bar{k}}}$; namely there exists the smallest number less than $m$ occurring in the sequence $(n_{i_k})_k$. By definition of $\al$ and the assumption that our sequence is descending, we have $n = n_{i_l}$ for all $l \geq \bar{k}$. Thus, $(i_l)_{l \geq \bar{k}}$ is an infinite and strictly increasing sequence of indices that witnesses $n$ to occur infinitely often in $(n_i)_i$ and this contradicts the hypothesis that RT$^1_{<\infty}$ fails; thus $\al$ is well-founded. The last step is to construct a bad sequence in $\al + \dots + \al$ whose elements have the form $(k,i)$ with $k \in \{0, \dots, m-1\}$ and $i \! \in \! \al$. Let us consider the sequence $\left( (m-1-n_i, i) \right)_i $ and let $i <_{\mN} j$ be two indexes; we claim that $(m-1-n_i,i) \nleq_{\times} (m-1-n_j, j)$ holds. If we have $ i >_{\al} j$, we are immediately done. Otherwise, if $i <_{\al} j$ holds, there are two cases: First, if we have $ n_i < n_j$, then this entails $m- 1 - n_i > m- 1 - n_j$ and, thus, our claim. In the latter case, if we have $n_i = n_j$ and $i >_{\mN} j$, then this immediately contradicts our assumption $i <_{\mN} j$. We can conclude that $\al + \dots + \al$ contains a bad sequence and is, therefore, not a well quasi-order.

\smallskip

\noindent ``\ita{3.} $\imp$ \ita{1.}'' We reason by contradiction. Fix $n$ and $Q$ such that $Q$ is {\wqo}, but $Q \, \dot{\cup} \dots  \dot{\cup} \, Q$ is not. Let $(\bar{q}_i)_i \in Q \, \dot{\cup} \dots  \dot{\cup} \, Q$ be an infinite bad sequence with $\bar{q}_i =(k_i, q_i)$, where $k_i \in \{0, \dots , n-1\}$ and $q_i \! \in \! Q$. Let us consider the function $f: \mN \frec \{0, \dots , n-1 \}$ defined as $f(i)=k_i$; by RT$^1_{<\infty}$ there exists $r \in \{0,\dots , n-1\}$ and an infinite increasing sequence $r_1 < r_2 < \dots $ such that $f(r_j)=r$ for all $j$. Since $(q_{r_j})_j$ is an infinite sequence in the well quasi-order $Q$, there exist two good indexes $s < t$ such that $q_{r_s} \leq_Q q_{r_t}$, for such indexes $\bar{q}_{r_s}=(r, q_{r_s}) \leq_{\dot{\cup}} (r, q_{r_t})=\bar{q}_{r_t}$; thus $Q \, \dot{\cup} \dots  \dot{\cup} \, Q$ is {\wqo}. \qed

\smallskip

For what concerns the ordinal analysis of WQP$(+^\omega)$, we have the two following results.

\begin{lem}\label{lem:+->+} \RCAs $\, \vdash \ \mbox{WQP}(+^\omega) \, \frec \, \mbox{WOP}(\g_{+})$.
\end{lem}

\proof Given a well-order $(\al, \leq_{\al})$, our goal is to prove that $\al \x \omega$ is a well-order; we actually prove that $\forall n\, \al \x n$ is a well-order, since this suffice over \RCAs. If $\al$ is a well-order, then $\al$ is a {\wqo} and, by WQP$(+^\omega)$, $\forall n\,$ $(n \x \al, \leq_{+})$ is a {\wqo} too. But \RCAs\, suffices to prove that $\leq_{+}$ is a linear order (given the linearity of $\leq_{\mN}$ and $\leq_{\al}$) and that $(n \x \al, \leq_{+})$ is order isomorphic to $\al \x n$, namely $(m_1, \al_1) \leq_{+} (m_2, \al_2)$ iff $(\al_1, m_1) \leq_{\al \x n} (\al_2, m_2)$ with respect to the well-order of the product $\al \x n$. \qed

\smallskip

Differently from WQP$(\times^\omega)$ and WOP$(\g_{\times})$, we do not have a proof for the separation between WQP$(+^\omega)$ and WOP$(\g_{+})$, nor for their equivalence.

\begin{prop}\label{prop:|+|} $|\mbox{\textsf{RCA}}_0^* + \mbox{WQP}(+^\omega)|=\omega^{\omega}$.
\end{prop}

\proof The proof is given by the following series of inequalities:
\[
\omega^{\omega} \leq |\mbox{\textsf{RCA}}_0^* + \mbox{WOP}(\g_{+})| \leq |\mbox{\textsf{RCA}}_0^* + \mbox{WQP}(+^\omega)| \leq |\mbox{\textsf{RCA}}_0| \leq \omega^{\omega}.
\]
The first inequalities derives from \Prop{prop:|g+|}, the second from \Lem{lem:+->+}, the third from the fact that \RCA\, proves both \RCAs\, and WQP$(+^\omega)$, the last one from the ordinal analysis of \RCA. \qed

Finally, we have the following.

\begin{prop}\label{prop:|U|} $|\mbox{\textsf{RCA}}_0 + \mbox{WQP}(\dot{\cup}^\omega)|=\omega^{\omega}$.
\end{prop}

\proof The proof is given by the following series of equalities:
\[
|\mbox{\textsf{RCA}}_0 + \mbox{WQP}(\dot{\cup}^\omega)| = |\mbox{\textsf{RCA}}_0 + \mbox{B}\Sigma^0_2| = \omega^{\omega},
\]
where the first equalities derives from \Theor{theor:4equiv}, and the second one from the case $n=2$ of \Cor{cor:Bsigma}. \qed

\smallskip

At the end of this section some general remarks on WOP's and WQP's are in order. First, we stress that, despite the use of the same symbols, ``$+$'' and ``$\times$'' have quite different meaning between well-orders, thus WOP, and {\wqo}, thus WQP. More precisely, in the context of well-orders, ``$+$'' and ``$\times$'' are binary operations on well-orders; whereas, talking about {\wqo}, ``$+$'' and ``$\times$'' are set operations that, joined with some quasi-orders construction, preserve {\wqo}. Secondly, whereas WOP's may have a uniform version and a ``standard version'', e.g. $\forall n\, \mbox{WOP}(\times^n)$ and WOP$(g_{\times})$, which can even be not equivalent (see \Cor{cor:WOPordinals}), this is not the case for WQP's; for example WQP$(\times^\omega)$ is defined simply as $\forall n\, \mbox{WQP}(\times^n)$.\footnote{In particular, although over $Q^\omega \putAs \{\mbox{infinite sequences on}\, Q \}$ a qo is easily definable, this is not in general a {\wqo} even when $Q$ it is; see Rado's counterexample \cite{Rado54}} Moreover, let us observe (cf. \Prop{prop:equiv} and \Prop{prop:KTln-WOP}) that well-ordering principles in the uniform version, i.e., in the form $\forall n\, \forall \fX \, [\dots]$, can be reduced to a standard form using $\sup_n$.

% 1. + e X diversi tra WOP e WQP 2. \forall n WOPn /WOPw versus \forall n WQPn /WQPw

In conclusion, the results of this section show that, differently from \ACA (cf. \Theor{theor:Equiv}), over \RCA, or weaker theories, {\wqo} and well-orders have quite different behaviours.

\section{Ordinal Analysis of $\Sigma_n$-Collection Schema}\label{sec:Bsigma}

This section, where we work mainly within first-order logic, is dedicated to the ordinal analysis of the \ita{Collection Schema} over a weak arithmetic theory. 

First of all, let us fix some notation. We adopt the standard language of first-order arithmetic, with symbols $0,1,+,\times,=,<$, possibly enriched with function symbols for a specific class of functions (in our case elementary functions) and with an extra uninterpreted function symbol $f$; $\Sigma_n$ and $\Pi_n$ formulas are defined analogously to second-order arithmetic, omitting the superscript $0$ or $1$ since no set variable is involved.

We take Elementary Recursive Arithmetic \ERA as base theory. \ERA is a first-order theory formulated in the standard language of first-order arithmetic extended with function symbols for all elementary functions; the set of \ita{elementary functions} is defined to be the smallest set of functions on the natural numbers containing the base functions $0, +, \times,  exp$, projections, and closed under the operations of composition and bounded recursion, see \cite{Rose85} for a detailed introduction. For what concerns its mathematical axioms, \ERA is axiomatized by: $(1)$ quantifier free defining equations for all elementary functions; $(2)$ usual axioms for $<$ and equality; $(3)$ the induction schema for quantifier free formulas. Moreover, bounded formulas in the language of \ERA are provably equivalent to quantifier free ones and it is known that \ERA itself admits a purely universal (or quantifier free) axiomatization.

For a set of formulas $\Gamma$, the related Induction Schema and Collection Schema, sometimes called Bounding Schema, are the following:
\[
\begin{array}{rl}\vspace{0.15 cm}

\mbox{I}\Gamma \putAs & \fhi(0) \et \forall k\, (\fhi(k) \frec \fhi(k+1) ) \, \frec \, \forall k\, \fhi(k) \ \mbox{for all}\ \fhi(x) \in \Gamma, \\

\mbox{B}\Gamma \putAs & \forall x \! \leq \! a\, \exists y\, \fhi(x,y) \, \frec \, \exists z\, \forall x \! \leq \! a\, \exists y\! \leq \! z\, \fhi(x,y) \ \mbox{for all}\ \fhi(x,y) \in \Gamma.

\end{array}
\]
We are particularly interested in the theories $ERA + \mbox{I}\Sigma_n$ and $ERA + \mbox{B}\Sigma_n$ which will be routinely denoted simply by $\mbox{I}\Sigma_n$ and $\mbox{B}\Sigma_n$; the distinction between the theory and the underlying axiom schema will be clear from the context. Furthermore, let us observe that the schema $\mbox{B}\Sigma_{n+1}$ is equivalent to $\mbox{B}\Pi_n$, and $\mbox{B}\Sigma_1$ is equivalent to the collection schema for quantifier free formulas.\footnote{For sake of completeness, note that $ERA + \mbox{B}\Sigma_1$ is essentially the first-order part of RCA$_0^*$, see \cite[Corollary 4.4]{SS:factorization}.}

The exploration of restricted collection principles in arithmetic began with Parsons \cite{Parsons:choice}, who used proof-theoretic approaches. This research was further advanced in the 1970s by Paris and his colleagues, after the pivotal role of collection principles in the model theory of arithmetic was recognized. In particular, Parsons established that the theories $\mbox{B}\Sigma_n$ and $\mbox{I}\Sigma_n$  form a strict hierarchy of fragments of Peano Arithmetic PA:
\[
ERA \, \subsetneq \, \mbox{B}\Sigma_1 \, \subsetneq \, \mbox{I}\Sigma_1 \, \subsetneq \, \mbox{B}\Sigma_2 \, \subsetneq \, \mbox{I}\Sigma_2 \, \subsetneq \mbox{B}\Sigma_3 \subsetneq \dots
\]
Paris \cite{Paris:hierarchy}, and independently H. Friedman, refined this picture proving that $\mbox{B}\Sigma_{n+1}$ is $\Pi_{n+2}$-conservative over $\mbox{I}\Sigma_n$, namely they proves the same $\Pi_{n+2}$ formulas. Given the relevance of this result in our analysis, we state it formally denoting with $\equiv_{\Gamma}$ the equivalence with respect to provability of formulas in the class $\Gamma$.

\begin{prop}\label{prop:conservativity} $\mbox{B}\Sigma_1 \equiv_{\Pi_2} ERA$ and, for all $n \geq 1$, $\mbox{B}\Sigma_{n+1} \equiv_{\Pi_{n+2}} \mbox{I}\Sigma_n$.
\end{prop}
Beklemishev \cite{Beklemishev:collection} not only proved this result using a syntactic approach that we exploit later, but he also further improved it calibrating the strength of $\mbox{B}\Sigma_{n+1}$ with respect of $\mbox{I}\Sigma_n$ in term of reflection principles.
 
As previously stated, our goal is to compute the $\Pi^1_1$ ordinal of the theory $\mbox{B}\Sigma_n$. In particular, let us define $\omega_1 \putAs \omega$ and $\omega_{n+1} \putAs \omega^{\omega_n}$, then we have the following.

\begin{theor}\label{theor:Bsigma} $|\mbox{B}\Sigma_1|=\omega^3$ and, for all $n \geq 2$, $|\mbox{B}\Sigma_n|=\omega_n$.
\end{theor}
The proof of this theorem, in particular for upper bounds, relies on the ordinal analysis of $\mbox{I}\Sigma_n(f)$ by Avigad in \cite{Avigad:without} and an appropriate extension of \Prop{prop:conservativity}.\footnote{The authors thank Fedor Pakhomov for suggesting this fruitful approach.}

For what concerns the former, Avigad studied fragments of first-order arithmetic, prominently \ERA and $\mbox{I}\Sigma_n$, enriched with an extra and uninterpreted function symbol $f$ which may appear in the definitions of other functions, e.g. using composition, but it is not required to satisfy any axiom. Syntactically, the classes $\Sigma_n(f)$ and $\Pi_n(f)$ can be defined as done for $\Sigma_n$ and $\Pi_n$, e.g. $\fhi$ is in $\Pi_2(f)$ if it has the form $\forall x\, \exists y\, A(x,y,f)$ with $A(x,y,f)$ having only bounded quantifiers. In the same vein, we can define the following extended versions of Induction and Collection schemas for $\Sigma_n(f)$ formulas:
\[
\begin{array}{rl}\vspace{0.15 cm}

\mbox{I}\Sigma_n(f) \putAs & \fhi(0) \et \forall k\, (\fhi(k) \frec \fhi(k+1) ) \, \frec \, \forall k\, \fhi(k) \ \mbox{for all}\ \fhi(x) \in \Sigma_n(f), \\

\mbox{B}\Sigma_n(f) \putAs & \forall x \! \leq \! a\, \exists y\, \fhi(x,y) \, \frec \, \exists z\, \forall x \! \leq \! a\, \exists y\! \leq \! z\, \fhi(x,y) \ \mbox{for all}\ \fhi(x,y) \in \Sigma_n(f).

\end{array}
\]
Let $ERA(f)$ be the theory sharing the same axioms of \ERA but with the syntax expanded with $f$. Similarly as before we denote by $\mbox{I}\Sigma_n(f)$ and $\mbox{B}\Sigma_n(f)$ the extensions of $ERA(f)$ with the schema $\mbox{I}\Sigma_n(f)$ and $\mbox{B}\Sigma_n(f)$, e.g. $\mbox{B}\Sigma_2(f)$ is the theory $ERA(f)$ enriched with the collection schema $\mbox{B}\Sigma_2(f)$.

Regarding ordinal analysis, Avigad calculated what here we call the ``computational ordinal'' of a given theory $T$ extending $ERA(f)$. More precisely, quoting \cite[$\S$4]{Avigad:without}, $\al$ is an upper bound of the computational ordinal of $T$ if there is a finitary proof of the following:
\begin{quote}
Whenever $T$ proves a $\Sigma_1(f)$ formula $\exists y\, \fhi(x, y, f)$, there is a $<\!\al$-recursive function $F$,\footnote{See \cite[$\S$3]{Avigad:without} for the precise definition of $<\!\al$-recursive function as well as for an appropriate encoding of countable ordinals within natural numbers.} such that for any $x, y$, and $f$, if $F(x, f) = y$, then $\fhi(x, y, f)$ is true.
\end{quote}
We rely on this specific type of ordinal analysis for two reasons: first, it implies an equal upper bound on the $\Pi^1_1$ ordinal we are interested in, i.e., if the computational ordinal of $T$ is almost $\al$, then $|T| \leq \al$; secondly, it involves the analysis exclusively of the $\Sigma_1(f)$ formulas of a theory, namely if $T$ and $T'$ are two theories extending $ERA(f)$ with the same set of $\Sigma_1(f)$ consequences, then they share the same computational ordinal. To convey a feeling on how such an ordinal analysis is pursued, and to clarify the latter of the two previous reasons, we restate, using our terminology and without proving it, \cite[Lemma 7.8]{Avigad:without} which is the key step to establish that $PRA(f)$ has as computational ordinal at most $\omega^{\omega}$:
\begin{quote}
Suppose $PRA(f)$ proves $\exists y\, \fhi(x,y,f)$ with $\fhi$ quantifier free; then there is a $<\! \omega^{\omega}$-recursive function $F(x,f)$ such that $ERA(f)$ proves $F(x,f)=y \frec \fhi(x,y,f)$.
\end{quote}
%For more complex theories, such as $\mbox{I}\Sigma_n(f)$, Avigad used diffe; nevertheless,
For what concerns our goal, we exploit the following estimation.
\begin{theor}\label{theor:Isigmaf} For $n$ greater than or equal to $1$, the computational ordinal of $\mbox{I}\Sigma_n(f)$ is at most $\omega_{n+1}$.
\end{theor}
\proof See \cite[Theorem 9.6]{Avigad:without}. \qed

To fruitfully apply this theorem to our case, namely $\mbox{B}\Sigma_n(f)$, we need to prove that $\mbox{B}\Sigma_{n+1}(f)$ and $\mbox{I}\Sigma_n(f)$ have, at least, the same $\Sigma_1(f)$ consequences. We actually prove something more.

\begin{lem}\label{lem:fconservativity} $\mbox{B}\Sigma_1 \! (f) \! \equiv_{\Pi_2(f)} \! ERA(f)$ and, for $n \!  \geq  \! 1$, $\mbox{B}\Sigma_{n+1}\!(f) \! \equiv_{\Pi_{n+2}(f)} \! \mbox{I}\Sigma_n(f)$.
\end{lem}
In order to derive this result, we adapt Beklemishev's proof \cite{Beklemishev:collection} checking that the correctness of the derivation is not disrupted by the addition of the new function variable $f$; Beklemishev's approach was chosen because, given its syntactic nature, is particularly suitable for this rewriting. For what concerns \Prop{prop:conservativity}%\footnote{As previously stated, \cite{Beklemishev:collection} contains also an analysis of $\mbox{B}\Sigma_{n}$ in terms of reflection principles.}
, Beklemishev's proof consists of three main steps:

\begin{enumerate}

 \item Instead of treating directly the collection schema $\mbox{B}\Sigma_n$, Beklemishev considers the following \ita{Collection Rule}:
\[
\Sigma_n\mbox{-CR}: \ \frac{\forall x\, \exists y\, \fhi(x,y)}{\forall x\, \exists y\, \forall u \! \leq \! x\, \exists v \! \leq \! y\, \fhi(u,v)}
\]
with $\fhi \in \Sigma_n$.

 \item He proves \cite[Theorem 3.1]{Beklemishev:collection} that the theory $\mbox{B}\Sigma_n$ is $\Pi_{n+1}$ conservative over $ERA + \Sigma_n\mbox{-CR}$. %in doing so, Beklemishev uses the equivalence between $\mbox{B}\Sigma_n$ and $\Pi_{n-1}$-FAC a finite form of axiom of choice (see later).
  
 \item Lastly, he checks \cite[Corollary 4.2]{Beklemishev:collection} that the theory $\mbox{I}\Sigma_n$ is equivalent to $ERA + \Sigma_{n+1}$-CR, i.e., they prove the same formulas. 

\end{enumerate}
Considering the extended collection rule $\Sigma_n(f)\mbox{-CR}\! : \ \frac{\forall x\, \exists y\, \fhi(x,y)}{\forall x\, \exists y\, \forall u \leq  x\, \exists v  \leq  y\, \fhi(u,v)}$ with $\fhi \! \in \! \Sigma_n(f)$, our adaptation of the previous proof consists in verifying that: i) the theory $\mbox{B}\Sigma_n(f)$ is $\Pi_{n+1}(f)$ conservative over $ERA(f) + \Sigma_n(f)\mbox{-CR}$, and ii) the theory $\mbox{I}\Sigma_n(f)$ is equivalent to $ERA(f) + \Sigma_{n+1}(f)$-CR. Since the aforementioned verification reduces to a long, often tedious, syntactic checking, for sake of brevity we treat explicitly only point i), proving it with \Prop{prop:Pequiv} after a series of preliminary results (which are themself an appropriate version of Beklemishev's preliminaries). The proof of point ii) is made in an analogous way, in particular adapting to the extended syntax the proof of \cite[Lemma 4.2]{Beklemishev:collection}.\footnote{Curiously, for proving Lemma 4.2 Beklemishev himself considers an extension of $ERA$ with unary function symbols.} 

\begin{prop}\label{prop:cont} $ERA(f) + \Sigma_{n+1}(f)\mbox{-CR}$ contains $\mbox{I}\Sigma_n(f)$.
\end{prop}

\proof Following \cite[Lemma 2.3]{Beklemishev:collection}, we derive an equivalent form of $\mbox{I}\Sigma_n(f)$, namely the \ita{strong $\Sigma_n(f)$ collection principle}:
\[
\mbox{S}\Sigma_n(f):\ \forall x\, \exists y\, \forall u \! \leq \! x \, ( \exists z \, \fhi(u,z,f) \, \frec \, \exists z \! \leq \! y \, \fhi(u,z,f)),
\]
where $\fhi \! \in \! \Sigma_n(f)$. The proof in \cite{HP98} of the equivalence between $\mbox{I}\Sigma_n$ and $\mbox{S}\Sigma_n$ can be easily extended to our enriched syntax.

To prove S$\Sigma_n(f)$ for a given formula $\fhi(x,y,f) \in \Sigma_n(f)$ using $\Sigma_{n+1}(f)\mbox{-CR}$, we proceed as follows. By appealing to excluded middle, clearly $ERA(f)$ proves
\[
\forall x\, \exists y\, (\fhi(x,y,f) \, \vel \, \forall z\, \neg \fhi(x,z,f)).
\]
Applying $\Sigma_{n+1}(f)\mbox{-CR}$, we derive
\[
 \forall x\, \exists y\, \forall u \! \leq \! x \, \exists v \! \leq \! y \, (\fhi(u,v,f) \, \vel \, \forall z\, \neg \fhi(u,v,f));
\]
thus, $ERA(f) + \Sigma_{n+1}(f)\mbox{-CR}$ proves
\[
\forall x\, \exists y\, \forall u \! \leq \! x \,  (\exists z\, \fhi(u,z,f) \, \frec \, \exists v \! \leq \! y\, \fhi(u,v,f)).  \tag*{\qed}
\]

Following Beklemishev's approach, it will be useful to adopt a sequent formulation for our Collection rule. We rely on a variant of Tait one-side calculus from \cite{Schwichtenberg77} in the rest of this section; in particular, sequents are finite sets of formulas (understood as disjunctions), negation is treated via de Morgan’s laws, etc. The next result enables the required sequent version of $\Sigma_n(f)\mbox{-CR}$.

\begin{prop}\label{prop:CR-schemas} The following rules are equivalent, in the sense that for any $R_1,R_2$ of them $ERA(f) + R_1 \equiv ERA(f) + R_2$.
\begin{enumerate}

\item 
\[
\Sigma_n(f)\mbox{-CR}: \ \frac{\forall x\, \exists y\, \fhi(x,y,f)}{\forall x\, \exists y\, \forall u \! \leq \! x\, \exists v \! \leq \! y\, \fhi(u,v,f)}
\]
for $\Sigma_n(f)$ formulas $\fhi(x,y,f)$.

\item 
\[
\Sigma_n(f)\mbox{-CR}': \ \frac{\forall x\, \exists y\, \fhi(x,y,f)}{\forall x\, \exists z\, \forall i \! \leq \! x\, \fhi(i,(z)_i,f)}
\]
for $\Sigma_n(f)$ or, equivalently, $\Pi_{n-1}(f)$ formulas $\fhi(x,y,f)$; $(z)_i$ denotes the $i$-th component of the finite sequence coded by $z$.

\item 
\[
\Sigma_n(f)\mbox{-CR}'': \ \frac{\Gamma, \forall x \! \leq \! t\, \exists y\, \fhi(x,y,f)}{\Gamma, \exists z\, \forall x \! \leq \! t\, \fhi(x,(z)_x,f)}
\]
where $\fhi(x,y,f) \! \in \! \Sigma_n(f)$ and $\Gamma$ is a set of $\Pi_{n+1}(f)$ formulas.

\end{enumerate}

\end{prop}

\proof \ita{1.} is easily derivable from \ita{2.} using the fact that $(z)_i \! \leq \! z$. For the other way round, assume that $ERA(f)$ proves $\forall x\, \exists y\, \fhi(x,y,f)$ with $\fhi \! \in \! \Pi_{n-1}(f)$. By $\Sigma_n(f)\mbox{-CR}$, we derive
\[
\exists z\, \forall i\! \leq \! x\, \exists y \! \leq \! z\, \fhi(i,y,f),
\]
thus,
\[
\exists z\, \forall i\! \leq \! x\, \exists y \! \leq \! (z)_i\, \fhi(i,y,f).
\]
By \Prop{prop:cont} $ERA + \Sigma_n(f)\mbox{-CR}$ contains I$\Sigma_{n-1}(f)$. By a theorem of Parsons \cite{Parsons:quantifier} (see also \cite{HP98}), which is extendable to our richer syntax, the least element principle for $\Delta_0(\Sigma_{n-1}(f))$ formulas is available in I$\Sigma_{n-1}(f)$; where by $\Delta_0(\Sigma_{n-1}(f))$ we mean formulas built from $\Sigma_{n-1}(f)$ ones using only boolean connectives and bounded quantifiers. Therefore, in $ERA + \Sigma_n(f)\mbox{-CR}$ there is a minimal $z$ such that
\[
\forall i\! \leq \! x\, \exists y \! \leq \! (z)_i\, \fhi(i,y,f).
\]
Since the standard coding of finite sequences is monotonic in each argument, $z$ is as required
\[
\forall i\! \leq \! x\, \fhi(i,(z)_i,f).
\]
\ita{2.} is derivable from \ita{3.} considering an empty $\Gamma$ and as $t$ a free variable. For the other direction, by dropping the outer universal quantifiers in $\Gamma$ (and introducing fresh free parameters) we may actually assume $\Gamma$ to be $\Sigma_n(f)$. Then the disjunction of $\Gamma$ can be moved inside the principal formula $\fhi$ of the choice rule. \qed

\begin{rem}\label{rem:B-FAC} Using almost the same strategy applied for the equivalence between 1. and 2. in \Prop{prop:CR-schemas}, it can be proved that, over $ERA(f)$, B$\Sigma_n(f)$ is equivalent to the following finite axiom of choice schema:
\[
\Pi_{n-1}(f)\mbox{-FAC}: \ \forall x \! \leq \! a\, \exists y \, \fhi(x,y,f) \ \frec \ \exists z\, \forall x \! \leq \! a\, \fhi\left(x,(z)_x,f\right),
\]
for $\Pi_{n-1}(f)$ formulas.
\end{rem}

Let us consider cut-free derivations of sequents of the form $\neg \mbox{B}\Sigma_n(f), \Pi(f)$ where $\Pi(f)$ denotes a set of $\Pi_{n+1}(f)$ formulas, and $\neg \mbox{B}\Sigma_n(f)$ is given by a finite set of negated instances of $\Sigma_n(f)$ collection schema, each of them having the form
\[
\exists z \exists u\, \left[ \forall x\! \leq \! u\, \exists y\, \fhi(x,y,z,f) \, \et \, \forall w\, \exists x \! \leq \! u\, \neg \fhi(x,(w)_x,z,f) \right]
\]
for some $\fhi \! \in \! \Pi_{n-1}(f)$; here we use the equivalence between B$\Sigma_n(f)$ and $\Pi_{n-1}(f)$-FAC ensures by \Rem{rem:B-FAC}. Moreover, let us denote the formula inside the square brackets with $B_{\fhi}(z,u,f)$.

Every sequent $\Gamma$ occurring in a cut-free proof of $\neg \mbox{B}\Sigma_n(f), \Pi(f)$ consists of formulas of the following two types: (i) formulas $\neg \mbox{B}\Sigma_n(f)$, or those of the form $\exists u\, B_{\fhi}(s,u,f)$ or $B_{\fhi}(s,t,f)$ for some $\fhi \! \in \! \Pi_{n-1}(f)$ and for some terms $s$ and $t$; (ii) subformulas of formulas from $\Pi(f)$, or proper subformulas of $B_{\fhi}(s,t,f)$. Finally, let $\Gamma^-$ denote the result of deleting all formulas of type (i) from such a sequent $\Gamma$, and notice that $\Gamma^-$ consists entirely of $\Pi_{n+1}$ formulas.

\begin{prop}\label{prop:cut} Let $\neg \mbox{B}\Sigma_n(f), \Pi(f)$ and $\Gamma$ be as before. If $\Gamma$ is cut-free derivable, then $\Gamma^-$ is derivable using $\Sigma_n(f)$-CR$''$ and first-order logic.
\end{prop}
\proof Let us proceed by induction on the height of the cut-free derivation $d$ of $\Gamma$. In particular, we only need to consider the case where a type (i) formula is introduced by the final inference in the derivation $d$; in this situation, we primarily need to handle formulas of the form $B_{\fhi}(s,t)$. In all other cases, once the $({\x})^-$-operation is applied, the premise and conclusion of the inference remain identical. Therefore, let us assume that the last step of derivation $d$ takes the following form:
\[
(\et) \ \ \frac{\forall x \! \leq \! t\, \exists y\, \fhi(x,y,s,f), \Delta \ \ \ \ \ \forall w\, \exists x \! \leq \! t\, \neg \fhi(x,(w)_x,s), \Delta}{B_{\fhi}(s,t), \Delta},
\]
with $\fhi \! \in \! \Pi_{n-1}(f)$.  Using the induction hypothesis, there are $\Sigma_n(f)$-CR$''$ derivations of
\begin{equation}\label{eq:1}
\forall x \! \leq \! t\, \exists y\, \fhi(x,y,s,f), \Delta^-,
\end{equation}
 and
\begin{equation}\label{eq:2}
\forall w\, \exists x\! \leq \! t\, \neg \fhi(x,(w)_x,s), \Delta^-.
\end{equation}
Since $\Delta^-$ consists of $\Pi_{n+1}(f)$ formula, applying $\Sigma_n(f)$-CR$''$ to \eqref{eq:1} we derive
\[
\exists w\, \forall x\! \leq \! t\, \fhi(x,(w)_x,s,f), \Delta^-.
\] 
The cut rule with sequent \eqref{eq:2} yields $\Delta^-$. \qed

\begin{prop}\label{prop:Pequiv} $\mbox{B}\Sigma_n(f) \, \equiv_{\Pi_{n+1}(f)} \, ERA(f) + \Sigma_n(f)\mbox{-CR}$.
\end{prop}

\proof Let $\fhi$ be a $\Pi_{n+1}(f)$ formula provable in $\mbox{B}\Sigma_n(f)$, and assume it is not provable in pure first-order logic (in which case it would be immediately provable also in $ERA(f) + \Sigma_n(f)\mbox{-CR}$). Then, there is a cut-free derivation of the sequent $\Gamma \putAs \neg \mbox{B}\Sigma_n(f), \Pi(f), \fhi$, where $\neg \mbox{B}\Sigma_n(f)$ is a finite set of instances of $\Sigma_n(f)$ collection, and $\Pi(f)$ is a finite set of instances of axioms of $ERA(f)$. Since $\Pi(f) \cup \{\fhi\}$ is a set of $\Pi_{n+1}(f)$ formulas, we can apply \Prop{prop:cut} to $\Gamma$ obtaining a derivation of $\Gamma^- = \Pi(f), \fhi$ in $\Sigma_n(f)\mbox{-CR}$ plus first-order logic (and thus in $ERA(f) + \Sigma_n(f)\mbox{-CR}$). Using repeatedly the cut rule and suitable instances of axioms of $ERA(f)$, we can eliminate all the formulas in $\Pi(f)$, obtaining $\fhi$. \qed

\smallskip

We can now prove our main result of this section.

\noindent \ita{Proof of \Theor{theor:Bsigma}:} For what concerns lower bounds, namely $\omega^3 \leq |\mbox{B}\Sigma_1|$ and $\omega_n \leq |\mbox{B}\Sigma_n|$, these derive immediately from the well-known ordinal analysis of \ERA and $\mbox{I}\Sigma_n$, together with the fact that $ERA \subset \mbox{B}\Sigma_1$ and $\mbox{I}\Sigma_n \subset \mbox{B}\Sigma_{n+1}$. Regarding upper bounds instead, we apply together \Theor{theor:Isigmaf} and \Lem{lem:fconservativity}. More precisely let $n \! \geq \! 2$, by \Lem{lem:fconservativity} $\mbox{B}\Sigma_n(f)$ and $\mbox{I}\Sigma_{n-1}(f)$ prove the same $\Pi_{n+1}(f)$, and thus also $\Sigma_1(f)$, formulas, and by Avigad's ordinal analysis they share the same computational ordinal. Therefore by \Theor{theor:Isigmaf}, $\mbox{B}\Sigma_n(f)$ has as computational ordinal at most $\omega_n$; since $\mbox{B}\Sigma_n \subset \mbox{B}\Sigma_n(f)$, we conclude $|\mbox{B}\Sigma_n| \leq \omega_n$. For $\mbox{B}\Sigma_1$ and $\mbox{B}\Sigma_1(f)$, the approach is completely analogous. \qed

The collection schema can be considered also in the context of second-order arithmetic in the following form:
\[
\mbox{B}\Sigma^0_n : \, \forall x \! \leq \! a\, \exists y\, \fhi(x,y) \, \frec \, \exists z\, \forall x \leq a\, \exists y \leq z\, \fhi(x,y),
\]
with $\fhi \in \Sigma^0_n$. From \Theor{theor:Bsigma}, we readily obtain the following corollary.

\begin{cor}\label{cor:Bsigma} For all $n \geq 2$, $|\mbox{\textsf{RCA}}_0 + B\Sigma^0_n|= \omega_n$.
\end{cor}
Extending previous results, in \cite{CSY:conservation} Chong et al. established that the principles %\footnote{Respectively the \ita{chain-antichain principle}, the \ita{ascending-descending principle}, and the \ita{cohesiveness principle}, see \cite{HS07} for their definitions.} 
CAC, ADS, and COH, coming from Ramsey and Computational theory (see \Subsec{subsec:preliminaries}), are all $\Pi^1_1$ conservative over \RCA$\, + \mbox{B}\Sigma^0_2$. Since theories with the same set of $\Pi^1_1$ consequences share the same $\Pi^1_1$ ordinal, we obtain:
\begin{cor}\label{cor:CAC} $|\mbox{\textsf{RCA}}_0 + \mbox{CAC}|=|\mbox{\textsf{RCA}}_0 + \mbox{ADS}|=|\mbox{\textsf{RCA}}_0 + \mbox{COH}|=\omega^{\omega}$.
\end{cor}

\section{Related and Future Work}\label{sec:works}

This paper is collocated at the crossroads between two well-established research topics: well-ordering principles and well quasi-orders. In this section, we summarize salient related papers and outline further developments.

\subsection{Related Literature}

Regarding well-ordering principles, the origin can be traced back to Girard \cite{Girard87} where the following equivalence was obtained:

\begin{theor}{(Girard\cite{Girard87})}\label{theor:ACAw} Over \RCA, \ACA\, is equivalent to WOP$(\la\fX.\omega^{\fX})$.\footnote{Here and below we use the $\la$-notation, namely $\la\fX.\phi(\fX)$, with $\phi(\fX)$ an ordinal term containing $\fX$, represents the function $\phi : \Omega \frec\Omega$ sending $\al$ to $\phi(\al)$.}
\end{theor}
Similar equivalences have been achieved for many other subsystems of second-order arithmetic.
\begin{theor}{(Afshari and Rathjen \cite{AR09}, Marcone and Montalbán \cite{MM11})} Over \RCA, \textsf{ACA}$_0^{+}$ is equivalent to WOP$(\la\fX. \ep_{\fX})$.
\end{theor}

\begin{theor}{(Friedman, Montalbán and Weiermann \cite{FMW07})} Over \RCA, \textsf{ATR}$_0$ is equivalent to WOP$(\la\fX.\fhi\fX0)$.
\end{theor}

\begin{theor}{(Rathjen \cite{Rathjen14})} Over \RCA, \textsf{ATR}$^{+}_0$ is equivalent to WOP$(\la\fX.\Gamma_{\fX})$.
\end{theor}

\begin{theor}{(Rathjen and  Vizcaíno \cite{RV15})} Over \RCA, \textsf{BI}$^{+}$ is equivalent to WOP$(\la\fX.\vartheta_{\fX})$.
\end{theor}

\begin{theor}{(Rathjen and  Thompson \cite{RT20})} Over \RCA, $\Pi^1_1$\textsf{-CA}$_0^{+}$ is equivalent to WOP$(\la\fX.\Omega_{\omega} \x \fX)$.
\end{theor}
For what concerns notation, in the previous results $T^+$ denotes, for a theory $T$ in the language of second-order arithmetic, $T$ extended with the following principle: ``every set is contained in a countable coded $\omega$-model of $T$''; see \cite{Arai20} or \cite{Rathjen22} for more details.

This connection between well-ordering principles and existence of suitable $\omega$-models has been explored by Michael Rathjen and co-authors in a series of papers \cite{AR09,Rathjen14,Rathjen22,RT20,RV15,RW11}; more recently, Anton Freund studied how $\Pi^1_1$-comprehension and well-ordering principles are related \cite{Freund18,Freund19,FR23}. The aforementioned results, e.g., $|\mbox{\textsf{ACA}}_0|=|\mbox{WOP}(\la\fX.\omega^{\fX})|=\ep_0$ \cite{Girard87}, $|\mbox{\textsf{ATR}}_0|=|\mbox{WOP}(\la\fX.\fhi\fX0)|=\Gamma_0$\cite{FMW07} on one side and $|\mbox{\textsf{ACA}}_0^{+}|=|\mbox{WOP}(\la\fX. \ep_{\fX})|$ \cite{MM11}, $|\mbox{\textsf{ATR}}_0^{+}|=|\mbox{WOP}(\la\fX.\Gamma_{\fX})|$ \cite{Rathjen14} on the other, suggest the existence of some general schemas. Such schemas have been exposed by Arai in \cite{Arai17,Arai20} where, for a normal function \g \ satisfying suitable term conditions \cite[Def.3 and Def.4]{Arai20} and its derivative $\g'$ (i.e., the function which enumerates the fixed points of $\g$), the following theorems have been proved.

\begin{theor}{(Arai\cite{Arai20})} $|\mbox{\textsf{ACA}}_0 + \mbox{WOP}(\g)|=\g'(0)=\min\{\al \,|\, \g(\al)=\al \}$.
\end{theor}

\begin{theor}{(Arai\cite{Arai20})} Over \ACA, the following are equivalent:

\begin{itemize}

\item WOP$(\g')$;

\item$(\mbox{WOP}(\g))^{+}$;

\end{itemize}
where $(\mbox{WOP}(\g))^{+}$ means that every set is contained in a countable coded $\omega$-model of \ACA$ \, + \,$WOP$(\g)$, see \cite[Definition 2]{Arai20} for a detailed definition.
\end{theor}
%Recently...

%\tbe \footnote{For the latter one of the standard references is \cite{Simpson09}.}

\medskip

Well quasi-orders have a much longer history, for which we refer to the now classic Kruskal's survey \cite{Kruskal72}. In a proof-theoretic perspective, among the extremely vast literature regarding {\wqo}, we draw our attention to the rich study of {\wqo} in constructive mathematics \cite{BBS24,BSB23,CF93,MR90,Powell20,RS93,SSW16,Veldman04} and, moreover, to the fruitful connection between {\wqo} and reverse mathematics. For example, over weak theories, such as \RCA, even the equivalence between different definitions of {\wqo} becomes non trivial as investigated by Marcone \cite{Marcone05,Marcone20} and Cholak et al. \cite{CMS04}. 
Also closure properties under {\wqo} play a relevant role in reverse mathematics with results such as the unprovability of Kruskal's theorem over \textsf{ATR}$_0$ \cite{RW93}, the unprovability of the Graph Minor theorem over $\Pi^1_1$\textsf{-CA}$_0$ \cite{FRS87}, or the following extension of \Theor{theor:ACAw}:

\begin{theor}{(Girard\cite{Girard87}, Hirst \cite{Hirst94}, Simpson\cite{Simpson88basis})}\label{theor:ACAHigman} Over \RCA\, the following are equivalent:

\begin{enumerate}

\item \ACA;

\item WOP$(\la\fX.\omega^{\fX}): \  \forall \fX\, \left[ \mbox{WO}(\fX) \frec \mbox{WO}(\fX^{\omega}) \right]$;

\item Higman's lemma$\, : \forall Q\, [Q\, \mbox{\wqo} \frec Q^* \, \mbox{\wqo}] $.

\end{enumerate}

\end{theor}
Remarkably, it was the failure of a closure property (the so called Rado's counterexample \cite{Rado54}) to lead Nash-William to introduce \ita{better quasi-orders}~\cite{Nash-W65trees}. 

\smallskip

For what concerns specific results in {\wqo} theory, two achievements stand out: Higman's lemma\footnote{The renowned lemma regarding finite strings is actually the corollary of a much wider result, see \cite{BS:wqo} for recent survey, and \cite{BW:relations} for a proof-theoretical analysis.} \cite{Higman52} and Kruskal's theorem \cite{Kruskal60}, whose forerunner was a conjecture by Vazsonyi. Regarding Kruskal's theorem in particular, after a plethora of applications in computer science, such as term rewriting \cite{Dershowitz82,Dershowitz87,DO88}, and mathematical logic, prominently ordinal analysis \cite{Gallier91,Schmidt79,VdMRW14,VdMRW16} and proof theory \cite{Simpson85nonprovability,Smith85}, the most recent results, due to Anton Freund and co-authors \cite{Freund20gap,FRW22,FU23}, regard mainly the so-called ``uniform Kruskal theorem'' which consists of an extension of the classical tree result to general recursive data types; involving Girard's dilator theory \cite{Girard81,Girard85}, this theorem has a strong categorical flavour.

Finally, for a recent collection that illuminates the multifaceted nature of {\wqo}, we refer to \cite{SSW20}.

\subsection{Future Work}\label{subsec:work}

There are three main directions for further developments concerning the topics treated in this paper:
\begin{itemize}

\item ordinal analysis over weak theories;

\item further analysis of well quasi-orders closure properties;

\item more precise proof-theoretic measures of Kruskal's theorem with labels.

%\item closure properties of wqo over weak theories.

\end{itemize}
Although we accomplished the ordinal analysis of WQP$(\g_{\times})$ and WOP$(\g_{+})$, some issues regarding ordinal analysis over \RCA\, and \RCAs\, remain, for example:

\begin{enumerate}

\item extend, under some suitable conditions over $\g$, the result of Arai \cite{Arai20} to \RCA, or even \RCAs, in order to obtain
\[
 |\mbox{\textsf{RCA}}_0^{(*)} + \mbox{WOP}(\g)|=\g'(0);
\]

\item in \cite{CMR19} Carlucci, Mainardi and Rathjen proved the following
\[
|\mbox{RCA}_0 + \mbox{WO}(\sigma) |= \sigma^{\omega} \ \mbox{for}\ \sigma \x \omega= \sigma,
\]
currently we lack a corresponding version of such result for \RCAs.

\end{enumerate}
Regarding the well quasi-orders closure properties considered in this paper, the main open problems are:
\begin{enumerate}

\item the relation between WQP$(\times^\omega)$ and RT$^2_{<\infty}$ over \RCA;

\item the relation between WQP$(+^\omega)$ and WOP$(\g_{+})$ over \RCAs.

\end{enumerate}
For what concerns Kruskal's theorem, the proof-theoretic investigations by Rathjen and Weiermann \cite{RW93} still point the way. Besides the proof-ordinal calculation, $\mbox{\textsf{RCA}}_0 \vdash \mbox{KT}(\omega) \sse \mbox{WO}\left( \vartheta (\Omega^{\omega}) \right)$, they sized the proof-theoretic strength of Kruskal's theorem for unlabelled trees in terms of reflection principles.

\begin{theor}(Rathjen and Weiermann \cite{RW93}) \RCA$\, \vdash \, \mbox{KT}(\omega) \sse \Pi^1_1\mbox{-}RF\!N(\Pi^1_2\mbox{-}BI_0)$.

\end{theor}
$\Gamma\mbox{-}RF\!N(T)$ denotes the theory $T$ extended with the uniform reflection principle $\forall x\, (Pr_T(\lceil \fhi (\bar{x})\rceil) \frec \fhi(x))$ for all formulas $\fhi$ in $\Gamma$ having at most one free variable, and $\Pi^1_2\mbox{-}BI_0$ is the theory \ACA$\, + \Pi^1_2\mbox{-}BI$, with $\Pi^1_2\mbox{-}BI$ the bar induction schema for $\Pi^1_2$ formulas (see \cite{RW93} for details). We refer to \cite{Beklemishev97,Beklemishev03,Beklemishev05,Rathjen94} for further readings regarding reflection principles, and in particular to \cite{JS99} for the connections between bar induction and reflection.

Consequentially, the next step is to obtain a similar classification also for \KTlw \ and $\forall n\,$\KTln. During the preparation of this paper, the two following conjectures arose:

\begin{conj}\phantom{placeholder}

\begin{itemize}

\item \RCA$\, \vdash \, \mbox{KT}_{\ell}(\omega) \sse \Pi^1_2\mbox{-}\omega RFN(\Pi^1_2\mbox{-}BI_0 \upharpoonright \Pi^1_3)$ \ \ [by F. Pakhomov],

\item \RCA$\, \vdash \, \forall n\, \mbox{KT}_{\ell}(n) \sse \Pi^1_2\mbox{-}RFN(\Pi^1_2\mbox{-}BI_0)$ \ \ \ \ \ [by A. Freund].

\end{itemize}

\end{conj}

Further work in this direction will also be required to extend the present results to other embeddability relations between trees, such as Friedman's gap condition \cite{Freund20gap,Kriz89}. %finally, for a recent survey on modern perspectives in Proof Theory see \cite{APW23}.

%\tbe

\section*{Conclusions}\label{sec:conclusion}

This paper has treated in the framework of ordinal analysis and reverse mathematics a series of Well-Ordering Principles, WOP, and Well Quasi-orders Principles, WQP, exploiting their strong connections for a joint study.

An extension of a previous result regarding ordinal analysis of well-ordering principles has been achieved in order to compute the proof-theoretic ordinals of different WOP's. More precisely, we extended the calculation, due to Arai \cite{Arai20}, of the proof-ordinal of the well-ordering principle relative to a normal function \g \ to the case of a weakly increasing ordinal function satisfying a side condition, i.e., $\g'(0)$ must be an epsilon number. The aforementioned extension has been subsequentially used to calculate the $\Pi^1_1$ ordinals of the WOP corresponding to the ordinal functions $\g_{\omega}(\fX) = \vartheta(\Omega^{\omega} \x \fX)$ and $\g_{\forall}(\fX) = \sup_n \vartheta(\Omega^n \x \fX)$. Using reverse mathematics techniques, we obtain also the $\Pi^1_1$ ordinals for the WOP relative to the functions $\g_{\times}(\fX)=\fX^{\omega}$ and $\g_{+}(\fX)=\fX \x \omega$.

These results allowed us to compute the proof-theoretic ordinals of \KTlw\ and $\forall n\,$\KTln, respectively the standard Kruskal's theorem for labelled trees and the conjunction of all the ``finite''versions for bounded branching trees; whereas the ordinal for the unlabelled case was already well known \cite{RW93}. In the same vein, the Well Quasi-orders closure Properties WQP$(+^\omega)$ and WQP$(\dot{\cup}^\omega)$, corresponding respectively to arbitrary sums and disjoint unions of {\wqo}, have been treated obtaining in both cases $\omega^{\omega}$ as proof-theoretic ordinal. Unexpectedly, this study unveiled connections with principles stemming from Ramsey and Computational theory, notably RT$^2_{<\infty}$ and~RT$^1_{<\infty}$.

Instrumentally to the ordinal analysis of WQP$(\dot{\cup}^\omega)$, and relying on previous achievements \cite{Avigad:without,Beklemishev:collection}, we compute also the proof-theoretic ordinal of the collection principles $\mbox{B}\Sigma_n$, seemingly missing in the literature.

Finally, our results have been situated within the current state of the art concerning well-ordering principles and well quasi-order closure properties, pointing toward further developments.

\subsection*{Acknowledgments}

The authors are deeply in debt to Toshiyasu Arai for the valuable suggestion regarding the extension of his result. The authors also thank Fedor Pakhomov and Anton Freund for their advices and for having pointed out the two reflection conjectures, as well as Patrick Uftring, Alberto Marcone and Lorenzo Carlucci for the fruitful discussions and the useful suggestion.% during the 2024 Logic Colloquium in Gothenburg.

The first author is member of the “Gruppo Nazionale per le Strutture Algebriche, Geometriche e le loro Applicazioni” (GNSAGA) of the ``Istituto Nazionale di Alta Matematica'' (INdAM).

\bibliographystyle{plain} %\small
\bibliography{{../Bibliografia/HigmanBiB}}   % name your BibTeX data base

\end{document}